\documentclass[preprint,12pt]{elsarticle}

\usepackage{graphicx}

\usepackage{amsmath}
\usepackage{amssymb}
\usepackage{soul}
\usepackage{url}
\usepackage{epstopdf}
\usepackage{rotating}
\usepackage{placeins}

\usepackage{lineno}		

\usepackage[caption=false,font=footnotesize]{subfig}
\usepackage{amssymb,amsfonts,mathrsfs,mathtools,amsthm}
\usepackage{booktabs}				
\usepackage{array}					
\usepackage{multirow}				
\usepackage{siunitx}					
\usepackage{flafter}					
\usepackage{comment}				
\usepackage[usenames,dvipsnames]{xcolor}	
\usepackage{flafter}					
\usepackage[margin=2.5cm]{geometry}	
\usepackage{enumitem}				
\usepackage{lipsum} 				
\usepackage{setspace}               

\usepackage{glossaries}	

\usepackage{nomencl}
\makenomenclature
\makeglossaries								

\usepackage{chngpage}

\usepackage{epstopdf}
\usepackage{lscape}
\usepackage{bm}
\usepackage{url}
\usepackage[font=footnotesize]{caption}
\usepackage{graphicx}
\usepackage{color}
\usepackage{algorithm,algorithmicx,algpseudocode}
\usepackage{amssymb}
\usepackage{nomencl}

\makenomenclature

\newacronym{cf}{CF}{capacity factor}
\newacronym{gw}{GW}{gigawatt}
\newacronym{gwh}{GWh}{gigawatt-hour}
\newacronym{hr}{HR}{heat rate}
\newacronym{ieee}{IEEE}{Institute for Electrical and Electronics Engineers}
\newacronym{mw}{MW}{megawatt}
\newacronym{mwh}{MWh}{megawatt-hour}

\newtheorem{assumption}{Assumption}
\newtheorem{remark}{Remark}

\usepackage{hyperref}	
\usepackage{cleveref}	

\usepackage{pifont}
\usepackage{graphicx}
\usepackage{threeparttable}
\usepackage{makecell}
\usepackage{footnote}

\usepackage{etoolbox}
\renewcommand\nomgroup[1]{%
  \item[\bfseries
  \ifstrequal{#1}{I}{Indices}{%
  \ifstrequal{#1}{V}{Variables}{%
  \ifstrequal{#1}{S}{Sets}{%
  \ifstrequal{#1}{P}{Parameters}{}}}}%
]}

\begin{document}

\begin{frontmatter}




\title{Incorporate Day-ahead Robustness and Real-time Incentives for Electricity Market Design\tnoteref{t1}}

\author[PSL]{Yi Guo}\ead{guo@eeh.ee.ethz.ch}

\author[PSL]{Xuejiao Han\corref{cor1}}\ead{xuhan@eeh.ee.ethz.ch}

\author[NREL]{Xinyang Zhou}\ead{xinyang.zhou@nrel.gov}

\author[PSL]{Gabriela Hug}\ead{hug@eeh.ee.ethz.ch}

\address[PSL]{Power Systems Laboratory, ETH Zürich, CH-8092 Zürich, Switzerland}
\address[NREL]{National Renewable Energy Laboratory, Golden, CO 80401, USA}

\cortext[cor1]{Corresponding author.}

\tnotetext[t1]{This work was supported by an ETH Postdoctoral Fellowship.}


\begin{abstract}
In this paper, we propose a two-stage electricity market framework to explore the participation of distributed energy resources (DERs) in a day-ahead (DA) market and a real-time (RT) market. The objective is to determine the optimal bidding strategies of the aggregated DERs in the DA market and generate online incentive signals for DER-owners to optimize the social-welfare taking into account network operational constraints. Distributionally robust optimization is used to explicitly incorporate data-based statistical information of renewable forecasts into the supply/demand decisions in the DA market.  We evaluate the conservativeness of bidding strategies distinguished by different risk aversion settings. In the RT market, a bi-level time-varying optimization problem is proposed to design the online incentive signals to tradeoff the RT imbalance penalty for distribution system operators (DSOs) and the costs of individual DER-owners. This enables tracking their optimal dispatch to provide fast balancing services, in the presence of time-varying network states while satisfying the voltage regulation requirement. Simulation results on both DA wholesale market and RT balancing market demonstrate the necessity of this two-stage design, and its robustness to uncertainties, the performance of convergence, the tracking ability and the feasibility of the resulting network operations.
\end{abstract}


\end{frontmatter}

\color{black}
\nomenclature[I]{\(t\)}{Index of time slots in the DA market}
\nomenclature[I]{\(k\)}{Index of time slots in the RT market}
\nomenclature[P]{\(T\)}{Number of time slots in the DA market}
\nomenclature[P]{\(\Delta T^{\textrm{DA}}\)}{Length of time slot in the DA market}
\nomenclature[P]{\(G_{t}\)}{Aggregated DA generation forecast at time $t$ in the DA market}
\nomenclature[P]{\(G^{\textrm{cap}}\)}{Generation capacity of DSO}

\nomenclature[P]{\(L_{t}\)}{Aggregated DA demand forecast at time $t$ in the DA market}
\nomenclature[V]{\(E^\text{BM+}_{t}/E^\text{BM-}_{t}\)}{Positive/negative imbalance quantities at time $t$}
\nomenclature[V]{\(\lambda^{\text{DA}}_{t}\)}{DA market clearing price for time $t$}
\nomenclature[V]{\(pr^{\text{BM+}}_{t}/pr^{\text{BM-}}_{t}\)}{Positive/negative imbalance prices at time $t$}
\color{red}\nomenclature[V]{\(E^{\text{DAs}}_{t}/E^{\text{DAb}}_{t}\)}{Bidding/offering quantities of DSO at time $t$}\color{black}
\nomenclature[P]{\(Tr^\text{max}\)}{Transmission capacity between the distribution and transmission grids}
\nomenclature[V]{\(\alpha^\text{DAs}_{t}/\alpha^\text{DAb}_{t}\)}{Offering/bidding price-quantity at time $t$ in the DA market}
\nomenclature[V]{\(E^\text{DA,D}_{t,m}\)}{Demand bidding quantity for block $m$ at time $t$}
\nomenclature[V]{\(E^\text{DA,O}_{t,j}\)}{Supply offering quantity for block $j$ at time $t$}
\nomenclature[V]{\(\lambda^\text{DA,D}_{t,m}\)}{Demand bidding price for block $m$ at time $t$}
\nomenclature[V]{\(\lambda^\text{DA,O}_{t,j}\)}{Supply offering price for block $j$ at time $t$}
\nomenclature[I]{\(j\)}{Index of bidding blocks in the DA market}
\nomenclature[P]{\(N^\text{j}\)}{Number of consumers' demand bidding blocks}
\nomenclature[I]{\(m\)}{Index of offering blocks in the DA market}
\nomenclature[P]{\(N^\text{m}\)}{Number of rival producers' offering blocks in the DA market}
\nomenclature[P]{\(E^{\text{DA,Omax}}\)}{Maximum DA offering quantity of the rival producers}
\nomenclature[P]{\(E^{\text{DA,Dmax}}\)}{Maximum bidding quantity of the consumers}

\nomenclature[V]{\(E^{\text{DAs,max}}_t/E^{\text{DAs,max}}_t\)}{Bidding/offering quantities of DSO at time $t$}
\nomenclature[V]{\(\bm{x}\)}{Compact vector collecting first-stage decisions}
\nomenclature[V]{\(\bm{y}\)}{Compact vector collecting second-stage decisions}
\nomenclature[V]{\(\mu_t^{\textrm{DAsmin}}/\mu_t^{\textrm{DAsmax}}\)}{Dual variable associated with lower/upper limits for DA dispatch supply quantity at time $t$}
\nomenclature[V]{\(\mu_t^{\textrm{DAbmin}}/\mu_t^{\textrm{DAbmax}}\)}{Dual variable associated with lower/upper limits for DA demand quantity at time $t$}

\nomenclature[V]{\(\mu_t^{\textrm{DA,Omin}}/\mu_t^{\textrm{DA, Omax}}\)}{Dual variable associated with lower/upper limits for DA offering quantity of the rival producers at time $t$}
\nomenclature[V]{\(\mu_t^{\textrm{DA, Dmin}}/\mu_t^{\textrm{DA, Dmax}}\)}{Dual variable associated with lower/upper limits for DA bidding quantity of the consumers at time $t$}

\nomenclature[P]{\(k^{0\sim3}, l^{0\sim3}\)}{Coefficients in linear decision rules}
\nomenclature[V]{\(u^{1,2}\)}{Auxiliary variables using linear decision rule.}
\nomenclature[V]{\(\tilde{y}(\bm{\delta},\bm{u})\)}{approximated resource decisions using linear decision rule}
\nomenclature[S]{\(\mathcal{N}_0\)}{Buses in the distribution network including the substation node}
\nomenclature[S]{\(\mathcal{N}\)}{Buses in the distribution network except the substation node}
\nomenclature[S]{\(\mathcal{E}\)}{Lines in the distribution network}

\nomenclature[V]{\(V_{i,k}\)}{Line-to-ground voltage at node $i$ at time $k$}
\nomenclature[V]{\(v_{i,k}\)}{Voltage magnitude at node $i$ at time $k$}
\nomenclature[V]{\(p_{i,k}/q_{i,k}\)}{Active/reactive power set-point of $i$-th DER at time $k$}
\nomenclature[P]{\(s_{i,k}^{\textrm{max}}\)}{Apparent power limit of $i$-th DER at time $k$}
\nomenclature[P]{\(p_{i,k}^{\textrm{min}}/p_{i,k}^{\textrm{max}}\)}{Lower/upper bounds of active power set-points of $i$-th DER at time $k$}
\nomenclature[P]{\(R/X\)}{Sensitivity matrices for power flow linearization}
\nomenclature[S]{\(\mathcal{X}_{i,k}\)}{Feasible set of $i$-th DER at time slot $k$ }
\nomenclature[V]{\(\alpha_{i,k},\beta_{i,k}\)}{Incentive signals for $i$-th DER at time $k$}
\nomenclature[P]{\(\Delta T^\textrm{RT}\)}{Length of time slot in the RT market}
\nomenclature[P]{\(E^\textrm{RT}_k\)}{RT balancing reference derived from the DA dispatch results at time $k$}
\nomenclature[P]{\(\underline{v}/\overline{v}\)}{Lower/upper bounds of voltage magnitude}
\nomenclature[P]{\(\tilde{v}\)}{Linearization coefficient of AC power flow}
\nomenclature[V]{\(\overline{\lambda}_k^{\textrm{RT}}/\underline{\lambda}_k^{\textrm{RT}}\)}{Dual variable associated with the upper/lower voltage limits at time $k$}


\printnomenclature[4cm]
\color{black}

\section{Introduction}
\label{sec:intro}
The continuing integration of distributed energy resources (DERs) in distribution networks, enhanced by the deployment of smart technologies at the end-user level, complicates balancing economic efficiency and system stability in distribution networks \cite{kroposki2020autonomous}. Such autonomous and intelligent DERs introduce both opportunities and challenges to the electricity market and electric power system operations. As the aggregations of DERs reach a substantial fraction of suppliers/consumers, they cannot be neglected as market participants in day-ahead (DA) and real-time (RT) markets any more. However, under current electricity market rules, DERs face high deliverable risks due to the unpredictable nature of renewable energy \cite{exizidis2019incentive,dvorkin2019chance,morales2013integrating}, which leads to security and reliability issues for distribution network operations. This motivates us to design a future electricity market mechanism that explicitly incorporates the stochasticity of aggregated DERs to manage these risks. We leverage distributionally robust DA bidding strategies and propose a fast incentive-based control for RT power balancing. These mechanisms account for the operational and economic objectives while also fulfilling constraints on voltages.

DERs are in general small-sized units that are connected to the distribution grid. Traditionally, end-consumers connected to the distribution grid face flat tariffs or two-tier tariffs (i.e., peak and off-peak tariffs). In this way, DER owners are exempt from additional costs in distribution system operation and maintenance resulting from DER injections or output forecast errors \cite{picciariello2015distributed}.
To promote better integration of DERs, attempts have been made to design local energy markets and new retail electricity tariff schemes. Local energy markets can be categorized into P2P energy markets and community-based markets \cite{zia2020microgrid}. Current retail electricity tariffs include flat tariff, time-of-use tariff and dynamic tariff, while the latter two time-based tariffs are proven to show few signs of cross-subsidization and better economic efficiency \cite{ansarin2020economic}. This work proposes a market framework for the optimal RT tariff design, while considering the tight connection between the DA wholesale market and local RT market with DERs. A detailed review along this line is provided in \cite{dutta2017literature}.

Existing works mainly focus on designing retail tariffs or pricing schemes for demand response programs.
A review on price-driven demand response programs is given in \cite{yan2018review}, which identifies that the price-signal can be an efficient tool for uncertainty and reliability management. In \cite{zhong2012coupon}, a coupon incentive-based demand response program is proposed on top of the flat retail electricity tariff.
The work in \cite{passey2017designing} designed a cost-reflective network tariff focusing on aligning the system's production and customers' demand peaks.

\begin{figure}
\centering
\includegraphics[width=5
in]{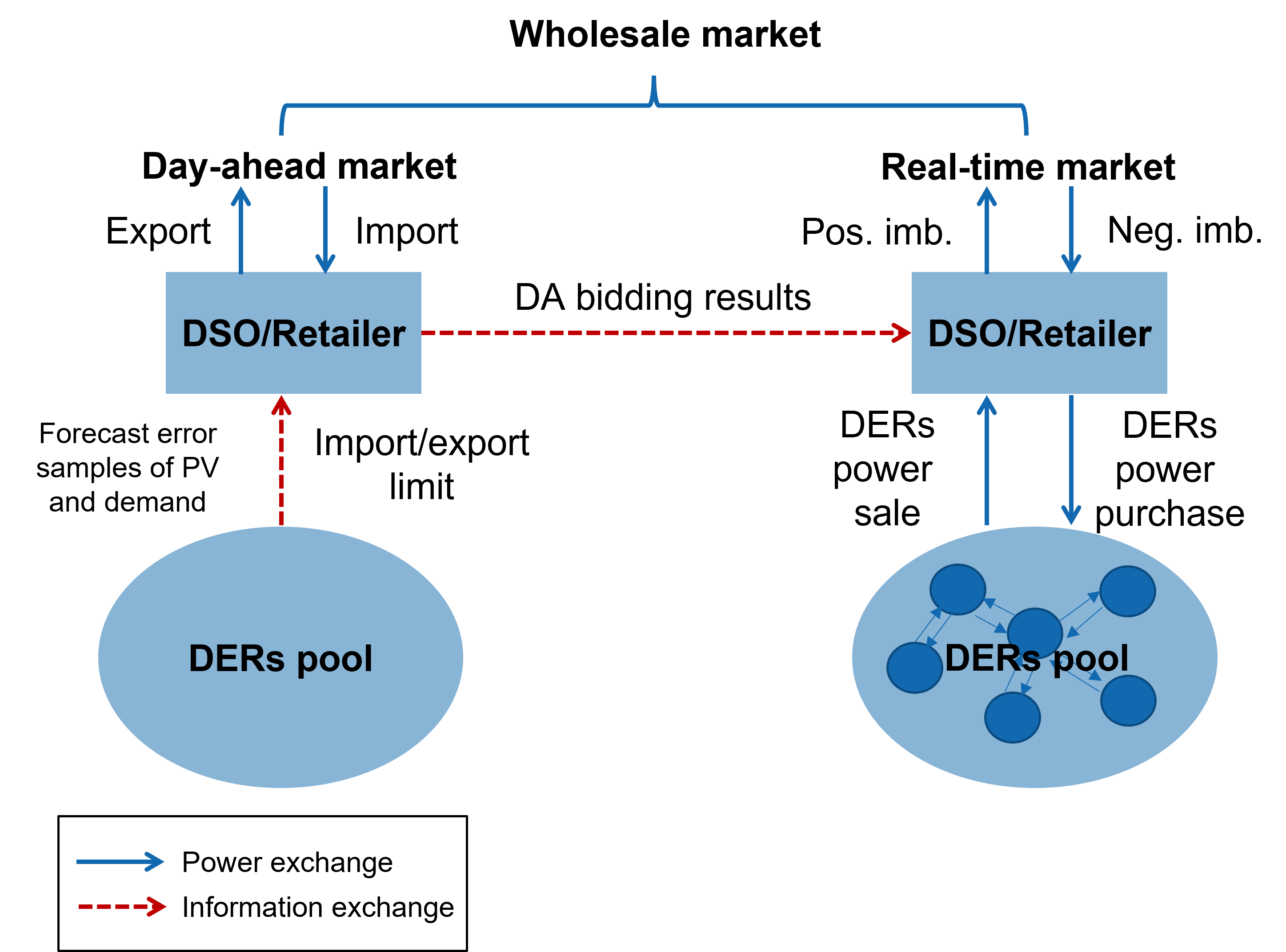}
\caption{The proposed two-stage market structure. In the DA market, the electricity supply and demand are scheduled based on the sampled forecast error dataset of renewable outputs. Considering the DA bidding strategies, the real-time incentive signals from the DSO enables to adjust the outputs of DER-owners such as to balance the power exchanges caused by the forecast errors. \color{black} The DER pool on the left indicates that all DERs in the distribution network act as a single market participant to bid into the DA wholesale market. The power exchanges between DERs are not considered at this stage. The connectivity between DERs on the right side of this diagram indicates that the DERs join the RT market with their own operational and economic objectives while taking into account the constraints imposed by the network.}
\label{fig_structure}
\end{figure}

While demand response programs are considered as flexibility providers that help to balance the system, the impacts brought by increasing DER penetrations can be either negative or positive, thus new pricing schemes that are based on a cost causation principle are required.
A review of network tariff design and incentives for DER owners can be found in \cite{picciariello2015distributed,eid2016managing}. 
A local market mechanism for a distribution network focusing on the external costs associated with voltage and line flow violations is presented in \cite{li2015market}. 
In \cite{guo2021online}, an online optimization framework that enables the  P2P market is introduced. A new business model for P2P energy sharing is proposed and comprehensively demonstrated in \cite{chen2021peer}.
Nevertheless, the aforementioned literature focuses on designing dynamic pricing schemes for the local market and ignores the connection to the wholesale market. The forecast errors of renewables will cause significant deviations of the RT dispatch from the DA generation/consumption schedule. The lack of the DA and RT markets' connection can reduce the available level of flexibility and lead to operations that violate network constraints. For a better interaction between the DA wholesale market and the RT local balancing market, we propose a two-stage electricity market consisting of a DA distributionally robust bidding process and an online distributed balancing algorithm. \color{black} The two-stage market mechanisms linking the DA market to the RT market for different types of DERs have been studied in the literature \cite{mirzaei2020novel,yang2017cvar,rahimiyan2015strategic,krishnamurthy2017energy,reddy2013optimal}. Compared to the existing works, the proposed electricity mechanism tackles the uncertainties from DERs by formulating the two-stage electricity market framework as a distributionally robust optimization problem for the slow time-scale of the DA market and an online optimization algorithm to cope with fast-changing DERs in the RT market. To the best of our knowledge, this is the first two-stage market mechanism that while linking the DA market with the RT market employs distinct and different algorithms for the different time scales. In addition, the cost of voltage regulation in the distribution network is also taken into account in the RT market decision, which is ignored in most of the literature. \color{black} The main contributions are as follows:

1) We formulate a two-stage electricity market problem for distribution networks with aggregated DERs. The framework is designed to enable the participation of DERs in the DA market in an aggregated way, and then uses a distributed incentive-based control strategy to enforce power and voltage constraints during RT operations.  In contrast to existing works, the proposed framework pursues the optimal power set-points of DERs in an online fashion while satisfying the network constraints and also accounts for the DA stochasticity realization by incorporating finite forecast samplings of renewable generations. The proposed overall market structure including the communication exchanges is presented in Fig.~\ref{fig_structure}. 
    
2) We formulate a stochastic DA bidding strategy using a bi-level model considering different levels of uncertainties and utilizing computationally tractable data-based stochastic optimization, i.e. distributionally robust optimization (DRO). Instead of assuming that the DERs' output forecasts follow prescribed probability distributions (e.g., Gaussian distribution), the proposed DRO market problem determines the optimal electricity supply/demand of a DER aggregator based on a forecast sampling dataset. These bidding strategies are robust to the \emph{worst-case} distribution within an ambiguity set, which consists of a group of probability distributions. This allows us to achieve superior out-of-sample performance of DA market results, efficiently avoiding overfitting the bidding to an available finite dataset.

3) A bi-level optimization is proposed to regulate local DERs for imbalance compensation in the RT market, in the presence of time-varying network conditions. The objective is to minimize the weighted sum of the imbalance costs for the DSO and the operational costs of DERs. The incentive-based signals for DER-owners are generated to adjust the local active/reactive set-points to balance the overall dispatch, while avoiding voltage constraint violations. An online implementation is proposed using a primal-dual gradient algorithm to achieve optimality from both the DSOs' and the DER-owners' perspectives. The effectiveness of the proposed market design is demonstrated on a wholesale market and an IEEE 37-node distribution network. 

The rest of the paper is organized as follow: Section \ref{sec:market_framework} introduces the DA wholesale market formulation using a distributionally robust optimization approach. Section \ref{subsec:marketstructure_B} describes the online incentive-based tariff design for a local RT balancing market. Section \ref{sec:num} provides the numerical results and Section \ref{sec:conclusions} concludes the paper.

\section{Day-Ahead Wholesale Market with Distributionally Robustness}\label{sec:market_framework}
Individual DER units are often small-sized and cannot participate in the wholesale market directly due to market restrictions such as minimum bidding quantity requirements. In this paper, we integrate multiple DER units within a distribution network into a single entity, i.e., aggregator, to bid in the wholesale market. 
\color{black}Note that only inverted-based distributed PV units are considered in this work, but it is straightforward to incorporate other distributed technologies, such as energy storage devices, demand response technologies or other distributed generators\footnote{\color{black}The proposed electricity market mechanism also allows to include DERs that lead to the bi-directional power flows as long as their models are linear, the cost functions are convex, and their feasibility sets are convex, closed and bounded. Besides, introducing energy storage devices into the proposed framework needs to include additional time-coupling constraints, but it does not change the property of optimality and convergence of the proposed algorithm.}.\color{black} We assume that the aggregator is namely the considered DSO and it bids into the DA wholesale market.
The other participants in the market, i.e. other aggregators, retailers, large scale power producers, etc., are modeled as demand and/or production bidders.
%

The objective of the DA optimization problem is to attain the DA dispatch decisions of the considered DSO, which is assumed to bid strategically into the market using a bi-level structure. All market participants other than the strategic DSO are assumed to be fully competitive and offer/bid their supplies/demands at their marginal prices. The bids of the considered DSO are optimized based on its forecasts of the aggregated market demand and supply curves. Solutions to this DA problem serve as the reference values for the local DERs represented by the DSO whose outputs are adjusted in the RT dispatch, which will be elaborated in Section~\ref{subsec:marketstructure_B}. 

To this end, we first divide each day equally into $T$ time slots indexed by $t\in\mathbb{T}=\{1,2,...,T\}$, where the length of each slot is $\Delta T^{\textrm{DA}}$. Given the aggregated DA generation and demand forecasts $G_{t} \in \mathbb{R}_{+}, L_{t} \in \mathbb{R}_{+}, \forall t\in\mathbb{T}$ of DERs, the considered DSO optimizes its DA bidding strategy  by solving a bi-level optimization problem. 
%
\subsection{Bidding Strategy Based on Bi-Level Optimization}
Assuming that the DSO bids strategically into the wholesale market. The interaction between the market and the DSO can be formulated as a leader-follower game, i.e., Stackelberg game \cite{zhang2009stackelberg}. We detail the problem formulation for the upper-level (UL) and lower-level (LL) problem for the leader and the follower in Sections \ref{subsec:UL} and \ref{subsec:LL}, respectively. 
\subsubsection{Upper-Level Optimization Problem}\label{subsec:UL}
The upper-level
(UL) optimization problem aims to minimize the DSO's cost
by optimizing its dispatch in the DA and balancing markets:
\begin{subequations}\label{eq:UL_problem}
\begin{eqnarray}
	\hspace{-5mm}& \underset{\tiny{\begin{subarray}{c}  E^\text{DAs/b,max}_{t},\\E^\text{BM+/-}_{t} \end{subarray}}}{\min} & \hspace{-5mm} \lambda^{\text{DA}}_{t}(E^{\text{DAb}}_{t}\!-\!E^{\text{DAs}}_{t})\!-\! pr^{\text{BM+}}_{t} E^{\text{BM+}}_{t}\!+\! pr^{\text{BM-}}_{t} E^{\text{BM-}}_{t}\!, \label{eq:UL_obj} \\ 
	\hspace{-5mm}&\text{s.t.} & E^{\text{DAs}}_{t} - E^{\text{DAb}}_{t} + E^{\text{BM-}}_{t} - E^{\text{BM+}}_{t} = G_{t} - L_{t}, \label{eq:UL_constraint0} \\ 
	\hspace{-5mm}&& 0\leq E^{\text{DAs,max}}_{t}\leq G^\text{cap},
	\label{eq:UL_constraint1} \\ 
	\hspace{-5mm}&& -Tr^\text{max} \leq E^{\text{DAs,max}}_{t} - E^{\text{DAb,max}}_{t} \leq Tr^\text{max},
	\label{eq:UL_constraint2} \\ 
    \hspace{-5mm}&& E^{\text{BM+}}_{t}, E^{\text{BM-}}_{t} \geq 0. \label{eq:UL_constraint-1} 
\end{eqnarray}
\end{subequations}
The objective function \eqref{eq:UL_obj} comprises the costs from the DA and the balancing markets, where $E^{\text{DAs}}_{t}$ (resp. $E^{\text{DAb}}_{t}$) represents the DA dispatched supply (resp. demand) quantity for time $t$, $E^{\text{BM+}}_t$ (resp. $E^{\text{BM-}}_t$) are the positive (resp. negative) imbalance quantities, and and $pr^{\text{BM+}}_t$ (resp. $pr^{\text{BM-}}_t$) are the positive (resp. negative) imbalance prices.
Constraint \eqref{eq:UL_constraint0} ensures the energy balance of the DSO. The bidding or offering quantities of the DSO, i.e., $E^{\text{DAs,max}}_t$ and $E^{\text{DAb,max}}_t$, are limited by the generation capacity $G^\text{cap}\in \mathbb{R}_{+}$ and the transmission capacity $Tr^\text{max}\in \mathbb{R}_{+}$ between the distribution and transmission grids in \eqref{eq:UL_constraint1}--\eqref{eq:UL_constraint2}.
%
\subsubsection{Lower-Level Optimization Problem}\label{subsec:LL}
The DA market clearing price $\lambda^{\text{DA}}_t$ in the UL objective function \eqref{eq:UL_obj} is the dual variable of the power balance of the lower-level (LL) market clearing problem that is solved by the market operator:
\begin{subequations}\label{eq:LL_problem}
\begin{eqnarray}
    &\underset{\tiny{\begin{subarray}{c}  \alpha^\text{DAs/b}_{t},E^\text{DAs/b}_t,\\ E^\text{DA,D}_{t,m},E^\text{DA,O}_{t,j} \end{subarray}}}{\max}\!& \alpha^\text{DAb}_{t} E^\text{DAb}_{t} + \sum^{N^\text{m}}_{m=1}\lambda^\text{DA,D}_{t,m}E^\text{DA,D}_{t,m} - \alpha^\text{DAs}_{t} E^\text{DAs}_{t} \nonumber \\ 
  && \qquad \qquad  - \sum^{N^\text{j}}_{j=1} \lambda^\text{DA,O}_{t,j} E^\text{DA,O}_{t,j},     \label{eq:LL_obj} \\ 
    &\text{s.t.}& \hspace{-5mm}
    \sum^{N^\text{m}}_{m=1}\hspace{-1mm}E^\text{DA,D}_{t,m}\hspace{-2mm}-\hspace{-1mm} \sum^{N^\text{j}}_{j=1}E^\text{DA,O}_{t,j}\hspace{-1mm}+\hspace{-1mm}E^{\text{DAb}}_{t} \hspace{-1mm}- \hspace{-1mm}E^{\text{DAs}}_{t}\!=\!0\!:\lambda^{\text{DA}}_t\hspace{-1mm}, \label{eq:LL_constraint0} \\ 
    && \hspace{-5mm} 0 \leq  E^\text{DAs}_{t}\leq E^\text{DAs,max}_{t}: \; \mu^\text{DAsmin}_{t}, \mu^\text{DAsmax}_{t}, \label{eq:LL_constraint1} \\ 
    && \hspace{-5mm} 0\leq  E^\text{DAb}_{t}\leq E^\text{DAb,max}_{t}: \; \mu^\text{DAbmin}_{t}, \mu^\text{DAbmax}_{t},\\ 
    && \hspace{-5mm} 0 \leq  E^\text{DA,O}_{t,j} \leq E^\text{DA,Omax}_{t,j}: \; \mu^\text{DA,Omin}_{t,j}, \mu^\text{DA,Omax}_{t,j}, \\ 
    && \hspace{-5mm} 0 \leq  E^\text{DA,D}_{t,m}\leq E^\text{DA,Dmax}_{t,m}: \; \mu^\text{DA,Dmin}_{t,m}, \mu^\text{DA,Dmax}_{t,m}. \label{eq:LL_constraint-1} 
\end{eqnarray}
\end{subequations}
For each time step, the aggregated offering and bidding curves consist of multiple supply offering and demand bidding blocks, i.e., offering and bidding price-quantity pairs on the step-wise offering and bidding curves. The objective of the LL problem is to maximize the total social welfare of the DA market, which is quantified by the demand bidding price-quantity pair $(\lambda^\text{DA,D}_{t,m},E^\text{DA,D}_{t,m}) \in \mathbb{R}^2_+$ for bidding block $m$ and the supply offering price-quantity pair $(\lambda^\text{DA,O}_{t,j},E^\text{DA,O}_{t,j}) \in \mathbb{R}^2_+$ for offering block $j$ of rival consumers and producers, respectively, and the offering (resp. bidding) price-quantity pair $(\alpha^\text{DAs}_{t}, E^\text{DAs}_{t}) \in \mathbb{R}^2_+$ (resp. $(\alpha^\text{DAb}_{t},E^\text{DAb}_{t}) \in \mathbb{R}^2_+$) of the DSO. Parameters $N^\text{m}$ and $N^\text{j}$ denote the number of consumers' demand bidding blocks and the number of rival producers' offering blocks, respectively.
Equality constraint (\ref{eq:LL_constraint0}) represents power balance constraint between the dispatched supply and demand quantities. Inequality constraints (\ref{eq:LL_constraint1})--(\ref{eq:LL_constraint-1}) enforce the minimum and the maximum limits of the offering and bidding quantities, where $E^{\text{DA,Omax}}\in \mathbb{R}_+^{N^\text{j}}$ and $E^{\text{DA,Dmax}} \in \mathbb{R}_+^{N^\text{m}}$ are denoted as the maximum DA offering quantity of the rival producers and the maximum bidding quantity of the consumers, respectively. Variables following colons after the equality and inequality constraints (\ref{eq:LL_constraint0})--(\ref{eq:LL_constraint-1}) are the corresponding dual variables.
\subsubsection{Combining LL and UL Problems}
As the bi-level optimization problem cannot be solved directly, reformulation is required. 
First of all, as the LL problem (\ref{eq:LL_problem}) is convex, the bi-level model can be formulated as a mathematical program with equilibrium constraints (MPEC) by replacing the LL problem the following set of Karush-Kuhn-Tucker (KKT) conditions:
\begin{subequations}
\begin{align}
& \alpha ^{\text{DAs}}_t - \lambda^{\text{DA}}_{t}+\mu^{\text{DAsmax}}_{t}-\mu^{\text{DAsmin}}_{t}=0 \; &\forall t, \label{eq:KKT_1} \\
& -\alpha ^{\text{DAb}}_t + \lambda^{\text{DA}}_{t}+\mu^{\text{DAbmax}}_{t}-\mu^{\text{DAbmin}}_{t}=0 \; &\forall t,  \label{eq:KKT_1b} \\
& \lambda^\text{DA,O}_{t,j}-\lambda ^{\text{DA}}_t+\mu^{\text{DA,Omax}}_{t,j}-\mu^{\text{DA,Omin}}_{t,j}=0 \; &\forall t,j,   \\
& -\lambda^\text{DA,D}_{tm}+\lambda ^{\text{DA}}_t+\mu^{\text{DA,Dmax}}_{t,m}-\mu^{\text{DA,Dmin}}_{t,m}=0 \; &\forall t,m,  \label{eq:KKT_3} \\
& 0\leq E^{\text{DAs}}_{t}\perp \mu^{\text{DAsmin}}_{t} \geq 0 \; &\forall t,   \label{eq:KKT_4} \\
& 0\leq E^{\text{DAb}}_{t}\perp \mu^{\text{DAbmin}}_{t} \geq 0 \; &\forall t,  \label{eq:KKT_4b}  \\
& 0\leq E^\text{DA,O}_{t,j}\perp \mu^{\text{DA,Omin}}_{t,j} \geq 0 \; &\forall t,j,  \\
& 0\leq E^\text{DA,D}_{t,m}\perp \mu^{\text{DA,Dmin}}_{t,m} \geq 0 \; &\forall t,m,   \\
& 0\leq (E^{\text{DAs,max}}_{t}-E^{\text{DAs}}_{t})\perp \mu^{\text{DAsmax}}_{t} \geq 0 \; &\forall t,  \label{eq:KKT_7} \\
& 0\leq (E^{\text{DAb,max}}_{t}-E^{\text{DAb}}_{t})\perp \mu^{\text{DAbmax}}_{t} \geq 0 \; &\forall t,  \label{eq:KKT_7b}\\
& 0\leq (E^{\text{DA,Omax}}_{t,j}-E^\text{DA,O}_{t,j})\perp \mu^{\text{DA,Omax}}_{t,j} \geq 0 \; &\forall t,j,   \\
& 0\leq (E^{\text{DA,Dmax}}_{t,m}-E^\text{DA,D}_{t,m})\perp \mu^{\text{DA,Dmax}}_{t,m} \geq 0 \; &\forall t,m, \label{eq:KKT_9}\\
& E^{\text{DAs}}_{t}-E^{\text{DAb}}_{t}+\sum^{N^\text{j}}_{j=1}E^\text{DA,O}_{t,j}-\sum^{N^\text{m}}_{m=1}E^\text{DA,D}_{t,m}=0 \;  &\forall t,j,m. \label{eq:KKT_10} \end{align}
\end{subequations}
where (\ref{eq:KKT_1})-(\ref{eq:KKT_3}) are stationary conditions. 
The resulting MPEC formulation is:
\begin{subequations}
	\begin{eqnarray}
	\hspace{-5mm}& \underset{\tiny{\begin{subarray}{c}  E^\text{DAs/b,max}_{t},\\E^\text{BM+/-}_{t} \end{subarray}}}{\min} & \hspace{-5mm} \lambda^{\text{DA}}_{t}(E^{\text{DAb}}_{t}\!-\!E^{\text{DAs}}_{t})\!-\! pr^{\text{BM+}}_{t} E^{\text{BM+}}_{t}\!+\! pr^{\text{BM-}}_{t} E^{\text{BM-}}_{t}\!,   \\ 
	\hspace{-5mm}&\text{s.t.} & \text{UL constraints~(\ref{eq:UL_constraint0})-(\ref{eq:UL_constraint-1})},\\
	\hspace{-5mm}&&\text{KKT conditions (\ref{eq:KKT_1})-(\ref{eq:KKT_10})}. 
	\end{eqnarray}
\end{subequations}

The resulting MPEC is non-linear due to 1) the complementarity conditions (\ref{eq:KKT_4})-(\ref{eq:KKT_9}) and 2) the term $\lambda^{\text{DA}}_{t}(E^{\text{DAb}}_{t}-E^{\text{DAs}}_{t})$ in the objective function. To convert the MPEC problem into a solvable MILP formulation, we first linearize the equations including the perpendicularity operator ''$\perp$'' using binary variables \cite{fortuny1981representation}:
\begin{subequations} 
\begin{align}
& 0 \leq E^{\text{DAs}}_{t} \leq M_1 u^{\text{DAsmin}}_{t} \quad &\forall t,  \\
& 0 \leq E^{\text{DAb}}_{t} \leq M_1 u^{\text{DAbmin}}_{t} \quad &\forall t,  \\
& 0 \leq E^\text{DA,O}_{t,j} \leq M_2 u^{\text{\text{DA,Omin}}}_{t,j} \quad &\forall t,j, \\
& 0 \leq E^\text{DA,D}_{t,m} \leq M_3 u^{\text{DA,Dmin}}_{t,m} \quad &\forall t,m, \\
& 0\leq \mu^{\text{DAsmin}}_{t} \leq M_4 (1-u^{\text{DAsmin}}_{t}) \qquad &\forall t,\\
& 0\leq \mu^{\text{DAbmin}}_{t} \leq M_4 (1-u^{\text{DAbmin}}_{t}) \qquad &\forall t,\\
& 0\leq \mu^{\text{DA,Omin}}_{t,j} \leq M_5 (1-u^{\text{\text{DA,Omin}}}_{t,j}) \qquad &\forall t,j,\\
& 0\leq \mu^{\text{DA,Dmin}}_{t,m} \leq M_6 (1-u^{\text{DA,Dmin}}_{t,m}) \qquad &\forall t,m,\\
& 0 \leq E^{\text{DAs,max}}_{t}-E^{\text{DAs}}_{t} \leq M_7 u^{\text{DAs,max}}_{t} \quad &\forall t, \\
& 0 \leq E^{\text{DAb,max}}_{t}-E^{\text{DAb}}_{t} \leq M_7 u^{\text{DAb,max}}_{t} \quad &\forall t, \\
& 0 \leq E^{\text{DA,Omax}}_{t,j}-E^\text{DA,O}_{t,j} \leq M_8 u^{{\text{DA,Omax}}}_{t,j} \qquad &\forall t,j, \\
& 0 \leq E^{\text{DA,Dmax}}_{t,m}-E^\text{DA,D}_{t,m} \leq M_9 u^{\text{DA,Dmax}}_{t,m} \quad &\forall t,m, \\
& 0\leq \mu^{\text{DAsmax}}_{t} \leq M_{10} (1-u^{\text{DAsmax}}_{t}) \qquad &\forall t,\\
& 0\leq \mu^{\text{DAbmax}}_{t} \leq M_{10} (1-u^{\text{DAbmax}}_{t}) \qquad &\forall t,\\
& 0\leq \mu^{\text{DA,Omax}}_{t,j} \leq M_{11} (1-u^{\text{DA,Omax}}_{t,j}) \quad &\forall t,j,\\
& 0\leq \mu^{\text{Dmax}}_{t,m} \leq M_{12} (1-u^{\text{DA,Dmax}}_{t,m}) \quad &\forall t,m,\\
& u^{\text{DAsmin}}_{t},u^{\text{DAbmin}}_{t},u^{\text{\text{DA,Omin}}}_{t,j},u^{\text{DA,Dmin}}_{t,m},u^{\text{DAsmax}}_{t},u^{\text{DAbmax}}_{t},u^{\text{DA,Omax}}_{t,j},u^{\text{DA,Dmax}}_{t,m} &\in \{0,1\}.
\end{align}
\end{subequations}
where $M_{1,2,...,12}$ are large enough constants.
Second, the non-linear term $\lambda^{\text{DA}}_{t}(E^{\text{DAb}}_{t}-E^{\text{DAs}}_{t})$ in the objective function is linearized by applying the strong duality theorem to the LL problem:
\begin{equation}\label{eq_DP}
\begin{aligned}
& \alpha^\text{DAb}_{t} E^\text{DAb}_{t} + \sum^{N^\text{m}}_{m=1}\lambda^\text{DA,D}_{t,m}E^\text{DA,D}_{t,m} - \alpha^\text{DAs}_{t} E^\text{DAs}_{t} - \sum^{N^\text{j}}_{j=1} \lambda^\text{DA,O}_{t,j} E^\text{DA,O}_{t,j}\\
= \; & \mu^{\text{DAsmax}}_{t} E^{\text{DAs,max}}_{t} + \mu^{\text{DAbmax}}_{t} E^{\text{DAb,max}}_{t} + \sum^{N^\text{j}}_{j=1}\mu^{\text{DA,Omax}}_{t,j} E^{\text{DA,Omax}}_{t,j} +  \sum^{N^\text{m}}_{m=1}\mu^{\text{DA,Dmax}}_{t,m} E^{\text{DA,Dmax}}_{t,m}
\end{aligned}
\end{equation}
By reformulating (\ref{eq:KKT_1})-(\ref{eq:KKT_1b}), (\ref{eq:KKT_4})-(\ref{eq:KKT_4b}) and (\ref{eq:KKT_7})-(\ref{eq:KKT_7b}), we obtain
\begin{subequations}
\begin{align}
&\alpha^{\text{DAs}}_{t} E^{\text{DAs}}_{t}= E^{\text{DAs}}_{t}(\lambda ^{\text{DA}}_t-\mu^{\text{DAsmax}}_{t}+\mu^{\text{DAsmin}}_{t}), \label{eq:temp1} \\
&\alpha^{\text{DAb}}_{t} E^{\text{DAb}}_{t}= E^{\text{DAb}}_{t}(\lambda ^{\text{DA}}_t+\mu^{\text{DAbmax}}_{t}-\mu^{\text{DAbmin}}_{t}), \label{eq:temp1b} \\
&E^{\text{DAs}}_{t} \mu^{\text{DAsmin}}_{t}=0,
\label{eq:temp2} \\
&E^{\text{DAb}}_{t} \mu^{\text{DAbmin}}_{t}=0,
\label{eq:temp2b} \\
&E^{\text{DAs}}_{t} \mu^{\text{DAsmax}}_{t}=E^{\text{DAs,max}}_{t} \mu^{\text{DAsmax}}_{t}, \label{eq:temp3} \\
&E^{\text{DAs}}_{t} \mu^{\text{DAbmax}}_{t}=E^{\text{DAb,max}}_{t} \mu^{\text{DAbmax}}_{t}.
\label{eq:temp3b}
\end{align}
\end{subequations}
Substituting (\ref{eq:temp2})-(\ref{eq:temp3b}) into (\ref{eq:temp1}) and (\ref{eq:temp1b}) yields
\begin{equation}
\begin{aligned}
\alpha^{\text{DAs}}_{t} E^{\text{DAs}}_{t}= \lambda ^{\text{DA}}_t E^{\text{DAs}}_{t}- \mu^{\text{DAsmax}}_{t} E^{\text{DAs,max}}_{t}, \\
\alpha^{\text{DAb}}_{t} E^{\text{DAb}}_{t}= \lambda ^{\text{DA}}_t E^{\text{DAb}}_{t} + \mu^{\text{DAbmax}}_{t} E^{\text{DAb,max}}_{t}
\end{aligned}
\end{equation}
and with (\ref{eq_DP}), we have
\begin{equation}
\begin{aligned}
\lambda^{\text{DA}}_t (E^{\text{DAb}}_{t}-E^{\text{DAs}}_{t})=& \sum^{N^\text{j}}_{j=1} \lambda^{\text{DA,O}}_{j} E^{\text{DA,O}}_{j}+\sum^{N^\text{j}}_{j=1}\mu^{\text{DA,Omax}}_{j} E^{\text{DA,Omax}}_{j}\\& -\sum^{N^\text{m}}_{m=1}\lambda^\text{DA,D}_{m} E^{\text{DA,D}}_{m}+\sum^{N^\text{m}}_{m=1}\mu^{\text{DA,Dmax}}_{m} E^{\text{DA,Dmax}}_{m}
\end{aligned}
\end{equation}
Following the linearization and reformulation process presented above, the bi-level problems (\ref{eq:UL_problem})--(\ref{eq:LL_problem}) can be reformulated as an MILP problem as follows:
\begin{subequations}\label{eq:DA_MILP}
\begin{eqnarray}
	&\hspace{-5mm}\underset{\tiny{\begin{subarray}{c}  \alpha^\text{DAs/b}_{t},E^{\text{DAs/b}}_{t},E^{\text{DAs,max}}_{t},\\E^\text{DA,O}_{t,j},E^\text{DA,D}_{t,m} \end{subarray}}}{\min } & \hspace{-3mm}
	\sum^{N^\text{j}}_{j=1} \big(\lambda^\text{DA,O}_{t,j} E^\text{DA,O}_{t,j}+\mu^\text{DA,Omax}_{t,j}E^\text{DA,Omax}_{t,j}\big)\nonumber\\[-3pt]
	&&\hspace{-5mm} 
	- \sum^{N^\text{m}}_{m=1}\big(\lambda^\text{DA,D}_{t,m}E^\text{DA,D}_{t,m}-\mu^\text{DA,Dmax}_{t,m}E^\text{DA,Dmax}_{t,m}\big) \nonumber \\
	&&\hspace{-5mm} 
	- pr^{\text{BM+}}_{t} E^{\text{BM+}}_{t}+ pr^{\text{BM-}}_{t} E^{\text{BM-}}_{t}, \label{eq:DA_MILP_obj} \\ 
	 &\hspace{-5mm}\text{s.t.}&
	 \hspace{-8mm}\text{UL constraints~(\ref{eq:UL_constraint0})--(\ref{eq:UL_constraint-1})} \label{eq:MILP_contraint0}, \\
   &&\hspace{-8mm}\text{Linearized reformulations of problem \eqref{eq:LL_problem}}. \label{eq:MILP_contraint-1} 
\end{eqnarray}
\end{subequations}
This combined equivalence of the bi-level problem can be directly solved using off-the-shelf commercial optimization solvers.
\subsection{Uncertainty Modelling}
In reality, the DA load and generation forecasts $L_t$ and $G_t$ in (\ref{eq:UL_constraint0}) are random variables.
Stochastic optimization (SO) \cite{pandvzic2013offering} and robust optimization (RO) \cite{liang2016robust,rahimiyan2016strategic} are among the most popular uncertainty modelling methods applied to optimize bidding strategies. In this paper we however use DRO to handle the uncertainty. This is because SO requires the knowledge of the specific uncertainty distribution and its computational complexity increases with the number of scenarios; although RO is often computationally tractable as it optimizes the decision considering the \textit{worst-case} scenario, performance of robust optimization is restricted by its conservativeness.

DRO was first developed for solving a single-product newsvendor problem considering a demand distribution characterized by its mean and variance in 1958 \cite{Scarf1958}. 
The method became popular again in recent years as it acts as an intermediary between SO and RO and achieves an acceptable trade-off between the optimality and the computational effort.
The distributionally robust DA bidding optimization presented here ensures that the bidding strategy is subject to the \textit{worst-case} distribution of generation/load forecast uncertainties within the ambiguity set $\mathcal{H}$. The ambiguity set collects a group of probability distributions of load and generation forecast errors $\delta_t \in \mathbb{R}$ and can be described by the following constraints
\begin{align}\label{eq:ambiguityset}:
\mathcal{H} = \left\{
\mathbf{H}: \begin{array}{l}
\mathbb{E}_{\mathbf{H}}[\delta_t]=0\\
\mathbb{E}_{\mathbf{H}}[|\delta_t|]\leq \zeta^1_{t}\\
\mathbb{E}_{\mathbf{H}}[(\delta_t)^2]\leq \zeta^2_{t}\\
Pr{(\delta_t\in \{\delta^\text{min}_t\leq \delta_t\leq \delta^\text{max}_t \}} ) =1\\
\end{array}
\right \},
\end{align}
where the first line ensures that the expectation of $\delta_t$ is zero. The second and third lines guarantee that the expected absolute deviation and the variance of $\delta_t$ are capped by $\zeta^1_{t}$ and $\zeta^2_{t}$, respectively. The last line limits all realizations of $\delta_t$ using the lower bound $\delta^\text{min}_t$ and the upper bound $\delta^\text{max}_t$. 

Following the principle of DRO \cite{Scarf1958}, we reformulate (\ref{eq:DA_MILP}) into a DA distirbutionally robust stochastic market problem:
\begin{subequations}\label{eq:minmax_DA_problem}
\begin{align}
& \inf_{\bm{x} } \sup_{\mathbf{H}\in\mathcal{H}} &&  \mathbb{E}^{\mathbf{H}}_{\bm{\delta}} \;
 [\; \Theta(\bm{x})+\phi(\bm{x,\delta}) \; ], \\
& \quad \textrm{s.t.} &&   \bm{Ax} + \bm{By}(\bm{\delta}) \leq \bm{D}(\bm{\delta}) \label{eq:minmax_DA_compactconstraint}, \\ 
& && \bm{\delta} \sim \mathbf{H}\in\mathcal{H},
\end{align} 
\end{subequations}
%
where terms $\Theta(\bm{x})$ and $\phi(\bm{x,\delta})$ correspond to the first-stage related part {\small$\sum^{N^\text{j}}_{j=1}  (\lambda^\text{DA,O}_{t,j} E^\text{DA,O}_{t,j}+\mu^\text{DA,Omax}_{t,j}E^\text{DA,Omax}_{t,j} ) - \sum^{N^\text{m}}_{m=1} (\lambda^\text{DA,D}_{t,m}E^\text{DA,D}_{t,m}-\mu^\text{DA,Dmax}_{t,m}E^\text{DA,Dmax}_{t,m} )$} and the second-stage related part {\small$pr^{\text{BM-}}_{t} E^{\text{BM-}}_{t}- pr^{\text{BM+}}_{t} E^{\text{BM+}}_{t}$} in the objective function \eqref{eq:DA_MILP_obj}, respectively. Constraint (\ref{eq:minmax_DA_compactconstraint}) is equivalent to (\ref{eq:MILP_contraint0})-(\ref{eq:MILP_contraint-1}).
The compact vectors $\bm{x}$ and $\bm{y}$ represent the first-stage and second-stage decisions (i.e., recourse decisions) defined by
\begin{equation*}
    \begin{aligned}
        & x:= [E^\text{DAs}_{t},E^\text{DAb}_{t},E^\text{DAs,max}_{t},E^\text{DAb,max}_{t},\alpha^\text{DAs}_{t},\alpha^\text{DAb}_{t}, \\
        & ~~~~~~~~\lambda^\text{DA}_{t},E^\text{DA,O}_{t,j},E^\text{DA,D}_{t,m}],\\
        & y: =[E^\text{BM-}_{t},E^\text{BM+}_{t}].
    \end{aligned}
\end{equation*}
%
%

\noindent Note that the ``min-max'' problem given in (\ref{eq:minmax_DA_problem}) can be reformulated as a minimization problem by taking the duality of the inner maximization problem \cite{Shapiro2001On}, however, the problem in general is still intractable as it requires solving recourse problems over all possible realizations of the uncertainty parameter $\bm{\delta}$ \cite{ben2004adjustable}. Thus, the concept of linear decision rule (LDR) \cite{goh2010distributionally} is applied to approximate the recourse decisions using an affine function of $\bm{\delta}$, i.e.,
\begin{subequations}\label{eq:LDR}
\begin{align}
	\tilde{E}^\text{BM+}_{t}&=k^0_{t}+k^1_{t}\delta_t+k^2_{t}u^1_{t}+k^3_{t}u^2_{t}, \\ 
	\tilde{E}^\text{BM-}_{t}&=l^0_{t}+l^1_{t}\delta_t+l^2_{t}u^1_{t}+l^3_{t}u^2_{t}, 
\end{align}
\end{subequations}
where $k^{0\sim3}_t \in \mathbb{R}$ and $l^{0\sim3}_t \in \mathbb{R}$ are coefficients to be optimized. The auxiliary variables $u^{1,2}_t \in \mathbb{R}$ are introduced to enhance the flexibility of the linear decision rule and guarantee the tractability of the problem \cite{wiesemann2014distributionally}, which are subject to the following additional constraints:
\begin{subequations}\label{eq:u_limits}
\begin{align}
& |\delta_t|\leq u^1_{t} \leq u^\text{1,max}_{t} = \max\{\delta^\text{max}_t,-\delta^\text{min}_t\}, \\
& (\delta_t)^2\leq u^2_{t} \leq u^\text{2,max}_{t} = \max\{(\delta^\text{max}_t)^2,(\delta^\text{min}_t)^2\}. 
\end{align} 
\end{subequations}

\noindent Substituting (\ref{eq:LDR}) into (\ref{eq:minmax_DA_problem}), we have:
\begin{subequations}\label{eq:minmax_DA_tractable}
\begin{align}
& \inf_{k^{0\sim3}_t,l^{0\sim3}_t} \sup_{\mathbf{H}\in\mathcal{H}} && \mathbb{E}^{\mathbf{H}}_{\bm{\delta}}\; \Theta(\bm{x})+\phi(\bm{x,\delta}), \\ 
& \textrm{s.t.} &&  \bm{x}\in\bm{X_f}, \\ 
& && \bm{Ax} + \bm{B\tilde{y}}(\bm{\delta},\bm{u})\leq \bm{D}(\bm{\delta}), \\ 
& && \bm{\delta} \sim \mathbf{H}\in\mathcal{H}, 
\end{align} 
\end{subequations}
where $\tilde{y}(\bm{\delta},\bm{u})$ denotes the approximated resource decisions using LDR. Eventually, problem
(\ref{eq:minmax_DA_tractable}) is tractable and it minimizes the DSOs' cost by optimizing the coefficients $k^{0\sim3}_t$ and $l^{0\sim3}_t$. 
As $\bm{\delta}$ and $\bm{u}$ are subject to constraints in the ambiguity set $\mathcal{H}$ as defined in \eqref{eq:ambiguityset} as well as \eqref{eq:u_limits}, problem \eqref{eq:minmax_DA_tractable} including the DRO constraints can be solved using the duality theory afterwards.

The outputs of the DA market, i.e., the dispatched DA bidding quantities of the DSO $\{E^\text{DAs}_t,E^\text{DAb}_t\}$, serve as inputs for the incentive-based RT balancing market. Details of the RT balancing market mechanism will be described in the following section. Note that the DA decisions are usually for an hourly resolution whereas actions taken based on the RT market outcomes are in the seconds to minutes time range to enable a timely tracking of time-varying loads and renewables. To resolve this inconsistency in temporal resolutions, we equally divide the DA bidding quantities of the DSO $\{E_t^\text{DAs},E_t^\text{DAb}\}$ into small portions to fit the fast balancing tasks, i.e. the power reference is assumed constant across all time slots in the RT problem that fall within one DA time slot. To avoid confusion with the notation, we use $k$ instead of $t$ as the time index when formulating the RT problem.
\section{Incentive-Based Real-Time Balancing Market}
\label{subsec:marketstructure_B}
In this section, we investigate a RT balancing market in distribution networks wherein both the DSO and DERs pursue their own operational and economic objectives. Again, the DSO corresponds to the entity that bids into the DA market as a representative of the local DERs. In the RT market, the DSO determines the optimal reward/payment of local DERs to encourage/discourage their network injection, such that the RT imbalance between the DA bidding quantity of the DSO and the RT 
output is minimized. A bi-level time-varying Stackelberg game-based optimization problem \cite{zhang2009stackelberg} is formulated to design the optimal incentive signals as well as the optimal set-points of the DERs. The controllability of the DERs are the operational set-points of the active and the reactive power. An online distributed algorithm is proposed to enable a computationally-efficient implementation.

\subsection{System Model}
Consider a distribution network operated by a DSO, denoted by a directed and connected graph $\mathcal{G}(\mathcal{N}_0,\mathcal{E})$, where $\mathcal{N}_0:= \mathcal{N}\cup\{0\}$ is the set of all ``buses" or ``nodes" with substation node 0 and $\mathcal{N}:= \{1,\dots,N\}$. The set $\mathcal{E} \subset \mathcal{N}\times\mathcal{N}$ includes ``links" or ``lines" for all $(i,j) \in \mathcal{E}$. Let $V_{i,k} \in\mathbb{C}$ denote the line-to-ground voltage at node $i\in\mathcal{N}$ at time $k$, where the voltage magnitude is given by $v_{i,k}:=|V_{i,k}|$. The set $\Omega$ includes all local DERs in a distribution network. Denote $p_{i,k}\in\mathbb{R}$ and $q_{i,k}\in\mathbb{R}$ as active and reactive power injections of DER at node $i$, respectively,  for all $i\in\Omega$ at time $k>0$. We denote $\mathcal{X}_{i,k}$ as the feasible set of active and reactive power $p_{i,k}$ and $q_{i,k}$ at node $i\in\Omega$ for all $k>0$. The set of operating set-points of DERs at node $i\in\Omega$ represents a convex envelop defined by
\begin{equation}\nonumber
    \mathcal{X}_{i,k}:=\Big\{(p_{i,k},q_{i,k}) : p_{i,k}^{\textrm{min}} \leq p_{i,k} \leq p_{i,k}^{\textrm{max}}, p_{i,k}^2 + q_{i,k}^2 \leq (s_{i,k}^{\textrm{max}})^2 \Big\},
\end{equation}
where $s_{i,k}^{\textrm{max}}$ is the apparent power limit of the DER at node $i\in\Omega$ at time $k$. Let $p_{i,k}^{\textrm{min}}$ and $p_{i,k}^{\textrm{max}}$ denote the lower and upper bounds of active power set-points of DER at node $i\in\Omega$ at time $k$. For PV inverter-based DERs, the feasible set $\mathcal{X}_{i,k}$ is constructed by the solar energy availability. For other devices, such as energy storage systems, small-scale diesel generators and variable frequency drives, the constraints can be altered to include their physical capacity limits in $\mathcal{X}_{i,k}$. We assume that the sets $\mathcal{X}_{i,k}$ are convex, closed and bounded for all $i\in\Omega$ for times $k\geq 0$. For future development, we define $\mathcal{X}_k:=\mathcal{X}_{1,k}\times\ldots\times\mathcal{X}_{N_{\Omega},k}$, where $N_{\Omega}$ denotes the cardinality of set $\Omega$.

To ensure that the optimal dispatch decisions of DERs are always feasible with respect to voltage constraints, we include the fundamental power flow equations in the RT balancing market design for distribution networks.
The AC power flow equations render the RT market problem nonconvex and NP-hard; in addition, they hinder the development of a computationally-affordable implementation. Here we instead use a linearization of the nonlinear AC power flow, which is given by
\begin{equation}
\label{eq:linear_powerflow}
    v_k = Rp_k + Xq_k + \tilde{v},
\end{equation}
where ${p}_k:=[p_{1,k},\ldots,p_{N,k}]^\intercal\in\mathbb{R}^N$ and ${q}_k:=[q_{1,k},\ldots,q_{N,k}]^\intercal\in\mathbb{R}^N$. The linearization parameters ${R}\in\mathbb{R}^{N\times N}$, ${X}\in\mathbb{R}^{N\times N}$ and $\bar{v}\in\mathbb{R}^{N}$ can be attained from various linearization methods, e.g., \cite{bolognani2015fast,baran1989network,guggilam2016scalable,christakou2013efficient,bernstein2017linear,dhople2015linear} and correspond to sensitivity matrices. 

\subsection{Real-time Incentive-based Market Problem}
The goals of our proposed incentive-based balancing market are 1) to explicitly take into account the inherent trade-offs between the renewable energy forecast errors in the DA dispatch results and the RT tariffs design for local DERs; 2) to coordinate DERs such as to fulfill the operational constraints (i.e., balancing and voltage regulations). Accordingly, two objectives are considered here to account for the different objectives for 
DERs and system operators.

\subsubsection{Costs for DERs}
The objective function for DERs at node $i\in\Omega$ comprises of the operational cost and the incentive cost, $J_{i,k}(p_{i,k},q_{i,k}) = J^{\textrm{Cost}}_{i,k}({p}_{i,k},{q}_{i,k}) + J^{\textrm{Inct}}_{i,k}({p}_{i,k},q_{i,k})$. The operational cost function $J^{\textrm{Cost}}_{i,k} \in \mathbb{R}_{+}$ is assumed to be quadratic and therefore convex, and can capture several objectives including ramping costs, small-scale thermal generation costs, active power losses and curtailment penalties. The incentive cost $J^{\textrm{Inct}}_{i,k} \in \mathbb{R}$ is a function of the incentive signals from the DSO to quantify the payment $J^{\textrm{Inct}}_{i,k}>0$ or reward $J^{\textrm{Inct}}_{i,k}\leq0$ for aggregated power injections of local DERs. We define the incentive costs for all DERs to be a linear affine function of the power dispatches, i.e., $J^{\textrm{Inct}}_{i,k}:= \alpha_{i,k} p_{i,k} + \beta_{i,k} q_{i,k}$. Both RT tariffs (incentives) $\alpha_{i,k} \in \mathbb{R}$ and $\beta_{i,k} \in \mathbb{R}$ and set-points of DERs $\{p_{i,k},q_{i,k}\}$ are decision variables in the following bi-level optimization problem. By this design, the incentives can be optimally adjusted by the DSO over time to continuously guide the power injections of DERs. Intuitively, the RT tariffs impact how much DERs should contribute to the overall RT dispatch of the DSO given the current network states and supply/demand conditions. For each time step $k$, the time-varying optimization problem for the DER at node $i\in\Omega$ can be expressed by
\begin{equation}\label{eq:DER_objective}
    \min_{(p_{i,k},q_{i,k})\in\mathcal{X}_{i,k}} J^{\textrm{Cost}}_{i,k}(p_{i,k}, q_{i,k}) + \alpha_{i,k}p_{i,k} + \beta_{i,k}q_{i,k}.
\end{equation}
Note that the above optimization is a convex quadratic program that determines the optimal set-points of DERs with
given incentive signals $\{\alpha_{i,k},\beta_{i,k}\}$. 

\subsubsection{Costs of Imbalance}
As the objective function for the DSO to be minimized, we define the discrepancy between the DA bidding results and the RT aggregated dispatch at the $k$-th time slot, i.e.,
\begin{equation}\label{eq:DSO_objectives}
    D_k(p_k | E^\textrm{RT}_k) = \bigg\|\sum_{i\in\mathcal{N}}p_k^i \Delta T^\textrm{RT} - E^\textrm{RT}_k\bigg\|_2^2,
\end{equation}
where $\Delta T^\textrm{RT}$ denotes the length of time slot during the RT operation. Let $E^\textrm{RT}_k$ denote the RT balancing reference derived from the DA dispatch results.
\begin{remark}(Timescale Mismatch and Connection)
We pursue the RT balancing market by tracking the DA dispatch decisions as a reference. The adjustments of DERs are determined every $\Delta T^\textrm{RT}$, but the DA reference values are given for every $\Delta T^{\textrm{DA}}$. To solve this timescale mismatch and build an appropriate connection between these two stages, we uniformly distribute the balancing task defined by $\{E_t^{\textrm{DAb}},E_t^{\textrm{DAs}}\}$ over the time slots in the RT operation. Given DA result at time $t$, the RT reference in \eqref{eq:DSO_objectives} therefore is
 \begin{equation}\label{eq:balancing_obj}
E_k^{\textrm{RT}} = \frac{E_t^{\textrm{DAs}} - E_t^{\textrm{DAb}}}{(\Delta T^{\textrm{DA}}/\Delta T^{\textrm{RT}})}, \quad t\cdot \Delta T^{\textrm{DA}} \leq k\cdot \Delta T^{\textrm{RT}} < (t+1)\cdot \Delta T^{\textrm{DA}},
\end{equation}
where the time index $t$ is reserved for the DA stage and index $k$ represents the time-step during the RT operation.
\end{remark}
\noindent Note that the social-welfare objective \eqref{eq:balancing_obj} indirectly connects the RT optimal balancing market to the quality of DA forecast. The stochastic modelling in the DA stage with different standard deviation settings possibly leads to overly conservative or riskier DA dispatch decisions. Our specific two-stage design allows DSOs to balance at a fast timescale to explicitly compensate the mismatch between the DA decisions and RT dispatches.

We introduce the following bi-level optimization problem, which captures both DER-oriented \eqref{eq:DER_objective} and network-oriented \eqref{eq:DSO_objectives} objectives
\begin{subequations}\label{eq:online_optimization_problem}
\begin{align}
    & \min_{\begin{subarray}{c} p_{k}, q_{k},\\ \alpha_{k}, \beta_{k}, v_{k} \end{subarray}} \sum_{i\in\Omega} J_{i,k}^{\textrm{Cost}}(p_{i,k},q_{i,k}) + \gamma D_k(p_{k}|E^\textrm{RT}_k),\\
    & \textrm{s.t.} \nonumber \\
    & (p_{i,k}, q_{i,k}) = \underset{p_{i,k}, q_{i,k}}{\arg\min}\  J^{\textrm{Cost}}_{i,k}(p_{i,k}, q_{i,k}) + \alpha_{i,k} p_{i,k} + \beta_{i,k} q_{i,k},  \label{eq:incentive_signal_design}\\
    & ~~~~~~~~~v_k = Rp_k + Xq_k + \tilde{v},\label{eq:power_flow_bi}\\
    & ~~~~~~~~~\underline{v} \leq v_k \leq \overline{v},\label{eq:voltage_bounds}\\
    & ~~~~~~~~~\forall i\in\Omega,\label{eq:incentive_signal design} 
\end{align}
\end{subequations}
where the constant $\gamma > 0$ is given based on the DA imbalance price. In the RT balancing market, the imbalance deviations happened on both sides (i.e., positive/negative imbalance quantities) will be penalized by the same cost, i.e., $\gamma$. The interaction between the imbalance costs and the cost functions of local DERs can also be seen as a Stackelberg game \cite{zhang2009stackelberg}.

The constraint \eqref{eq:incentive_signal_design} models an embedded optimization problem for the DER at node $i\in\Omega$ with given incentive signals $\alpha_{i,k}$ and $\beta_{i,k}$. The linear power flow in \eqref{eq:power_flow_bi} maps the active and reactive set-points $(p_{k},q_{k})$ to voltage magnitude with sensitivity matrix $(R, X)$ at any time slot $k$.  The vectors $\underline{v}\in\mathbb{R}^N$ and $\bar{v}\in\mathbb{R}^N$ represent the lower and upper limits of voltage magnitude. In the RT market, the DSO is responsible for power balancing and voltage regulation at any time. This time-varying problem is posed and solved every $\Delta T^{\textrm{RT}}$ for the ``best" incentive signals and optimal set-points for DERs', while at the same time regulating the voltage. It is a challenge to solve problem \eqref{eq:online_optimization_problem} in real time, not only because of the non-convex nature of the problem but also because it requires continued communication between DSO and DERs due to the time-varying situation (i.e., supply/demand variations). To tackle this issue, we leverage a gradient approach to approximate the solution of problem \eqref{eq:online_optimization_problem} in an online distributed fashion.
\subsection{Online Distributed Algorithm: A Gradient Approach}
Online gradient-based approaches deal with optimization problems that have incomplete or time-varying input information (parameters). The decisions are implemented over time without fully solving the optimization problem, and it aims to tradeoff optimality, communication effort and computational efficiency. Such algorithms have been applied and discussed for power systems and other applications in \cite{dall2016optimal,tang2017real,guo2021online,bernstein2015composable,gan2016online,bernstein2015design}. In this paper, we employ an online gradient-based approach to solve the proposed RT incentive-based balancing market problem taking into account the fast-changing renewable power output. While developing the online algorithm, we firstly start with a convex relaxation of the original problem \eqref{eq:online_optimization_problem}, such that the optimization problem of local DERs, i.e. \eqref{eq:incentive_signal_design}, is ignored. Then we show that a primal-dual gradient approach together with a specific design of the incentive signal updates can attain the optimum of the original problem \eqref{eq:online_optimization_problem}. Replacing the embedded constraint \eqref{eq:incentive_signal_design} by the operational feasible region, we obtain
\begin{subequations}\label{eq:online_relaxation}
\begin{align}
    & \min_{\begin{subarray}{c} p_{i}, q_{i}, v_{i} \end{subarray}} && \sum_{i\in\Omega}J_{i,k}^{\textrm{Cost}}(p_{i,k},q_{i,k}) + \gamma D_k(p_{k}|E^\textrm{RT}_k),\\
    & \textrm{s.t.} &&v_k = Rp_k + Xq_k + \tilde{v},\label{eq:power_flow}\\
    &  &&\underline{v} \leq v_k \leq \overline{v}: \underline{\lambda}_k^{\textrm{RT}},  \overline{\lambda}^{\textrm{RT}}_k,\label{eq:voltage_bounds}\\
    & &&(p_{i,k},q_{i,k}) \in \mathcal{X}_{i,t}, \forall i\in\mathcal{N},\label{eq:incentive_signal design} 
\end{align}
\end{subequations}
where $\underline{\lambda}_k^{\textrm{RT}}\in\mathbb{R}^N_{+}$ and $\overline{\lambda}_k^{\textrm{RT}}\in\mathbb{R}^N_{+}$ are the dual variables associated with the lower and upper voltage constraints, respectively. We make the following assumptions.

\begin{assumption}\label{ass:differentiability}
For any time $k>0$, the local objective functions of DERs, $J^{\textrm{Cost}}_{i,k}(p_{i,k},q_{i,k}), \forall i\in\Omega$ are continuous differentiable and strongly convex functions of $p_{i,k}$ and $q_{i,k}$, and their first-order derivatives are bounded within their operation regions, i.e., $\nabla J^{\textrm{Cost}}_{i,k}(p_{i,k},q_{i,k}) \leq M_{J}, \forall i\in\Omega$. The imbalance cost function $D_k(p_k |E^\textrm{RT}_k) = \|\sum_{i\in\mathcal{N}}p_{i,k} \Delta T^{\textrm{RT}} - E^\textrm{RT}_k\|_2^2$ is continuously differentiable, convex and with first-order derivative bounded by given DA bidding strategies $E^\textrm{RT}_k$, i.e., $\nabla D_k(p_k |E^\textrm{RT}_k)\leq M_{D}$.
\end{assumption}
\begin{assumption}[Slater's condition]\label{ass:slater}
For any time $k>0$, there exists a feasible point located within the operating region $(p_{k},q_{k})\in\mathcal{X}_{k}$, so that
    \begin{equation*}
        \underline{v} \leq Rp_k + Xq_k + \tilde{v} \leq \overline{v}.
    \end{equation*}
\end{assumption}

\begin{remark}[Optimal Condition]
Under Assumptions \ref{ass:differentiability} and \ref{ass:slater}, the solution of \eqref{eq:online_relaxation} along with the incentive signals $(\alpha_{k}^{*},\beta_{k}^*)$ defined by
\begin{subequations}\label{eq:incentive_signal_define}
\begin{align}
    \alpha_k^* & = R\left(\underline{\lambda}_k^{\textrm{RT},*} - \overline{\lambda}_k^{\textrm{RT},*} + \gamma \nabla_{p_k}D_k(p_k^*|E_k^{\textrm{RT}})\right), \label{eq:incentive_signal_define_alpha}\\
    \beta_k^* & = X\left(\underline{\lambda}_k^{\textrm{RT},*} - \overline{\lambda}_k^{\textrm{RT},*} \right),
\end{align}
\end{subequations}
is the global solution of the original problem \eqref{eq:online_optimization_problem}.
\end{remark}
\noindent The proof is omitted here, which follows a similar derivation as given in \cite{zhou2017incentive}.

We now develop an online gradient-based algorithm to solve the RT balancing market problem in \eqref{eq:online_optimization_problem}. Initially, consider a regularized Lagrangian function of the relaxed problem  \eqref{eq:online_relaxation} given by
\begin{equation}\label{eq:L_function}
\begin{aligned}
    & \mathcal{L}_k^\eta\left(p_k, q_k, \underline{\lambda}_k^{\textrm{RT}}, \overline{\lambda}_k^{\textrm{RT}}\right) \\
    & = \sum_{i\in\mathcal{N}}J_{i,k}^{\textrm{Cost}}(p_{i,k},q_{i,k}) + \gamma D_k(p_{k}|E^\textrm{RT}_k) + (\underline{\lambda}_k^{\textrm{RT}})^\intercal \left(\underline{v} - v_k\right)  \\
    & +  (\overline{\lambda}_k^{\textrm{RT}})^\top \left(v_k - \underline{v}\right) - \frac{\eta}{2}\left(\|\underline{\lambda}_k^{\textrm{RT}}\|_2^2 + \|\overline{\lambda}_k^{\textrm{RT}}\|_2^2\right),
\end{aligned}
\end{equation}
where a small positive constant $\eta>0$ is predefined. The Tikhonov regularization term $-\frac{\eta}{2} \left( \|\underline{\lambda}_k^{\textrm{RT}}\|_2^2 + \|\overline{\lambda}_k^{\textrm{RT}}\|_2^2 \right)$ facilitates the convergence performance. 
To solve \eqref{eq:online_relaxation} in an online fashion, we formulate the time-varying saddle-point problem  
\begin{equation}\label{eq:saddle-point}
    \max_{\underline{\lambda}_k^{\textrm{RT}}\in\mathbb{R}^N_{+},\overline{\lambda}_k\in\mathbb{R}^N_{+}}  \min_{(p_t,q_i)\in\mathcal{X}_{i,k}}  \quad \mathcal{L}_k^\eta\left(p_k, q_k, \underline{\lambda}_k^{\textrm{RT}}, \overline{\lambda}_k^{\textrm{RT}}\right).
\end{equation}
As $\eta$ is small, the primal-dual gradient-based approaches can be applied to \eqref{eq:saddle-point} to reach an approximate solution of the original problem \eqref{eq:online_relaxation} but with better convergence. The optimality discrepancy due to the regularization terms has been explicitly discussed in \cite{koshal2011multiuser}. Together with the incentive signals defined in \eqref{eq:incentive_signal_define}, we have the following online iterative updates at time $k$:
\begin{subequations}\label{eq:online_algorithm}
\begin{align}
    p_{k+1} & = \left[p_{k} - \epsilon_p \left(\nabla_{p_k} \sum_{i\in\mathcal{N}}J_{i,k}^{\textrm{Cost}}(p_{i,k},q_{i,k}) + \alpha_{k}\right) \right]_{\mathcal{X}_{k}}, \label{eq:online_update_p}\\
    q_{k+1} & = \left[q_{k} - \epsilon_q\left(\nabla_{q_k} \sum_{i\in\mathcal{N}}J_{i,k}^{\textrm{Cost}}(p_{i,k},q_{i,k})+\beta_k\right)  \right]_{\mathcal{X}_{k}},\label{eq:online_update_q}\\
    \underline{\lambda}_{k+1} & = \bigg[\underline{\lambda}_{k} + \epsilon_{\lambda}\left(\underline{v} - v_k - \eta\underline{\lambda}_k^{\textrm{RT}} \right) \bigg]_{\mathbb{R}_+},\label{eq:online_update_dual_1}\\
    \overline{\lambda}_{k+1} & = \bigg[\overline{\lambda}_{k} + \epsilon_{\lambda}\left(v_k - \overline{v} -  \eta\overline{\lambda}_k^{\textrm{RT}} \right) \bigg]_{\mathbb{R}_+},\label{eq:online_update_dual_2}\\
    \alpha_{k+1} & = R\bigg(\underline{\lambda}_{k+1}^{\textrm{RT}} - \overline{\lambda}_{k+1}^{\textrm{RT}} + \gamma \nabla_{p_k}D_k(p_k|E_k^{\textrm{RT}}) \bigg),\label{eq:online_update_alpha}\\
    \beta_{k+1} & = X\bigg(\underline{\lambda}_{k+1}^{\textrm{RT}} - \overline{\lambda}_{k+1}^{\textrm{RT}} \bigg),\label{eq:online_update_beta}\\
    v_{k+1} &~~~~ \textrm{updates based on sensor measurement}, \label{eq:online_v_feedback}
\end{align}
\end{subequations}
\noindent where $\epsilon_p$, $\epsilon_q$ and $\epsilon_{\lambda}$ are the positive constant step-sizes for primal and dual updates. The operator $[\cdot]_{\mathcal{X}_k}$ projects onto the feasible set $\mathcal{X}_k$. The operator $[\cdot]_{\mathbb{R}_+}$ projects onto the nonnegative orthant. The above iterations \eqref{eq:online_algorithm} are performed over time for time steps $k>0$ with time-varying updates of the problem formulation \eqref{eq:online_optimization_problem}. Due to the space limitations, we omit the discussions of the convergence performance and tracking capability but it will be demonstrated in the simulation section.
Here, we mostly focus on enhancing the connections and exploring the tradeoffs between the DA market results and the RT implementation. 

The updates \eqref{eq:online_algorithm} are of a distributed nature and therefore can be implemented in a distributed way. For any $k>0$, the DERs $i\in\Omega$ update their operational points $(p_{i,k},q_{i,k})$ locally through \eqref{eq:online_update_p}--\eqref{eq:online_update_q} based on their individual incentives $\{\alpha_{i,k},\beta_{i,k}\}$. Note that it is not necessary for DERs to broadcast their own cost functions $J_{i,k}^{\textrm{Cost}}$ and operational region $\mathcal{X}_{i,k}$. Similarly, the topology information of distribution networks $(R, X)$ and the DA bidding results $\{E^{\textrm{DAb}},E^{\textrm{DAs}}\}$ are not required for the local computations. The system operators require the set-points of local DERs and RT voltage measurements \eqref{eq:online_v_feedback} to update the dual variables \eqref{eq:online_update_dual_1}--\eqref{eq:online_update_dual_2} and generate the incentive tariffs \eqref{eq:online_update_alpha}--\eqref{eq:online_update_beta}. Hence, the DSO and local DERs  can coordinate effectively and on a fast timescale under this time-varying setup. This enables privacy preservation of DERs and DSO and a minimal communication load to achieve both network-oriented and DER-oriented objectives. Algorithm \ref{algorithm:real-time_market_design} below summarizes our proposed online incentive-based market algorithm. 

\begin{remark}(Interpretation of Incentive signals) The terms in the tariffs defined in \eqref{eq:incentive_signal_define} can be uniquely assigned to specific incentive goals. The terms that are functions of the dual variables associated with voltage limits, i.e.,
\begin{equation}\nonumber
    \alpha_{k}^{\textrm{V}} = R^{\top}\left(\underline{\lambda}_k^{\textrm{RT}} - \overline{\lambda}_k^{\textrm{RT}}\right), \beta_{k}^{\textrm{V}} = X^{\top}\left(\underline{\lambda}_k^{\textrm{RT}} - \overline{\lambda}_k^{\textrm{RT}}\right),
\end{equation}
incentivize local DERs to contribute to the voltage regulation. The last term of \eqref{eq:incentive_signal_define_alpha}, i.e.,
\begin{equation}\nonumber
    \alpha_{k}^{\textrm{DSO}} = \gamma R\nabla_pD_k(p_k|E_{k}^{RT}),
\end{equation}
quantifies how much the DSO encourages/discourages their DERs to adjust the set-points to contribute to the balancing objective. The weighted summation of these two parts leads to the final incentive information in \eqref{eq:incentive_signal_define} for which the DSO needs to tradeoff the network performance and the RT market response by defining the parameter $\gamma$.
\end{remark}

\begin{algorithm}
    \caption{(Online Incentive-based Market Implementation)}\label{algorithm:real-time_market_design}
    \begin{algorithmic}[1]
        \Require[S0] DA bidding decisions $\{E^{\textrm{DAb}},E^{\textrm{DAs}}\}$. Initialization of set-points of DERs $\{p_0,q_0\}$, incentive signals $\{\alpha_0,\beta_0\}$,  dual variables $\{\underline{\lambda}_0, \overline{\lambda}_0 \}$ and voltage profile $v_0$.
        \While{$k = 1:T^{\textrm{RT}}$}
            \State[S1] Network DSO performs the dual updates \eqref{eq:online_update_dual_1}-\eqref{eq:online_update_dual_2} for voltage regulation.
            \State[S2] Network DSO calculates the incentive signals \eqref{eq:online_update_dual_1}-\eqref{eq:online_update_dual_2} and pass them to local DERs.
            \State[S3] Local DERs $i\in\Omega$ perform updates of power set-points \eqref{eq:online_update_p}-\eqref{eq:online_update_q}
            \State[S4] Local DERs $i\in\Omega$ implement the power set-points $\{p_{k+1},q_{k+1}\}$.
            \State[S5] Network DSO collects the voltage magnitude measurement $v_{k+1}$.
        \EndWhile
\end{algorithmic}
\end{algorithm}

\begin{figure}[h]
    \centering
    \includegraphics[width=4.0in]{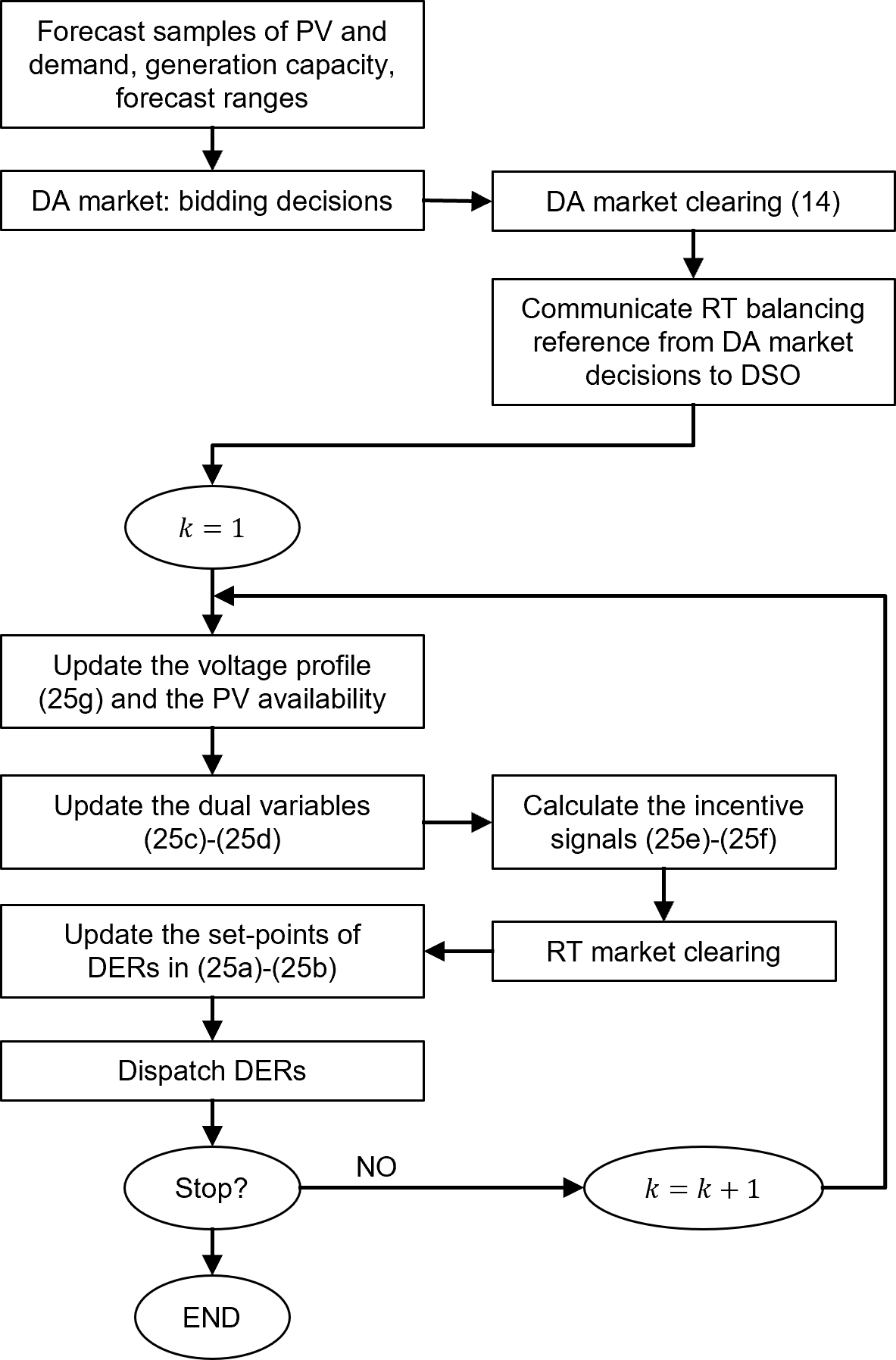}
    \caption{\color{black}The decision sequence of the proposed two-stage electricity market mechanism.}
    \label{fig:flowchart}
    \end{figure}

\color{black} Finally, the decision sequence of the proposed electricity market mechanism can be summarized in the following steps:
        \begin{itemize}
            \item[1)] The DSO attains the forecast samples of PV and demand and determines the available generation capacity and forecast ranges.
            \item[2)] The DSO determines and submits its hourly aggregated bidding decisions to the DA market for the next day.
            \item[3)] The market operator clears the DA market.
            \item[4)] The DSO receives the hourly RT balancing reference from the DA market decision.
            \item[5)] The DSO attains the PV availability and measures the voltage magnitude \eqref{eq:online_v_feedback} at time $k$.
            \item[6)] The DSO performs the dual update \eqref{eq:online_update_dual_1}-\eqref{eq:online_update_dual_2} at time $k$.
            \item[7)] The DSO calculates the incentive signals \eqref{eq:online_update_dual_1}-\eqref{eq:online_update_dual_2} based on the current voltage profile, PV availability, and the difference between the DA bidding results and the RT aggregated dispatch at time $k$.
            \item[8)] The DSO sends the incentive signals to DERs at time $k$.
            \item[9)] The DSO clears the RT market at time $k$.
            \item[10)] DERs update and implement the set-points of power injections \eqref{eq:online_update_p}-\eqref{eq:online_update_q} at time $k$.
        \end{itemize}
The above steps 5)-10) are repeated over time until the termination. The flowchart in Fig.~\ref{fig:flowchart} visualizes the above decision sequence. \color{black}

\section{Numerical Results}\label{sec:num}
In this section, we provide numerical results for the proposed algorithms. First, we study the DA market problem and discuss the results thereof before then demonstrating the workings of the RT balancing mechanism. 
\subsection{Day-ahead Market Data and Results}\label{subsec:DA}
The DA biddings of other market participants are assumed to be perfectly forecasted and are modeled as parameters using the DA market supply and demand curves based on the Nord Pool market clearing data of 2018.
For each hour the original supply (offering) and demand (bidding) curves consist of up to nearly 1000 blocks. Due to the computational burden, the original supply (offering) and demand (bidding) curves are approximated focusing on the bids and offers near the original market clearing point. The resulting approximated supply and demand curves for each hour consist of a maximum of 79 blocks.
The considered DSO is modeled as a new prosumer to the market. In other words, the DSO enters the market by adding its offers and bids to the existing offering and bidding curves.  
The total Nord Pool system-level bidding and offering quantities, which cover the area of several countries, are scaled to simulate the case for one transmission system.

Balancing market prices are modeled such as to guarantee that the DSO can only sell (purchase) electricity in the balancing market at a price lower (higher) than the DA market price:
\begin{equation}\label{eq:price_imb}
pr^\text{BM+}_{t} = a_1 \cdot (pr^\text{DA}_{t}-p_1),
pr^\text{BM-}_{t} = a_2 \cdot (pr^\text{DA}_{t}+p_2),
\end{equation}
where $pr^\text{DA}_t$ are the original DA market clearing prices. Constants $a_1$ and $a_2$ are set to 0.7 and 1.7, price adjustments $p_1$ and $p_2$ are set to 15 EUR/MWh and 20 EUR/MWh, respectively.  

The DA forecast error of the PV power output is modelled as a random variable. For each time step, 1000 samples of PV forecast errors $\delta_t$ are generated randomly under a Gaussian distribution with zero mean and standard deviations $\sigma$ equaling to 0, 10\% and 20\% of the PV generation capacity. We evaluate the DA bidding decisions under these three different standard deviations of PV forecasts. Note that as the generation outputs are non-negative values and are limited by the PV capacity, unrealistic forecasts that fall out of this range are adjusted to the respective bounds. After replacing the unrealistic forecast samples, the empirical distribution supported by the sampling dataset is no longer a Gaussian distribution. This motivates us to leverage the DRO to make market decisions based on a group of distributions to capture the real unknown data-generated distributions. We assume that the upper and lower bounds of $\delta$ used to construct the uncertainty set for DRO, the empirical mean, the mean absolute deviation and standard deviation used to build the ambiguity set for DRO are calculated based on the sampled data.

\begin{figure}
\centering
\includegraphics[width=0.7\textwidth]{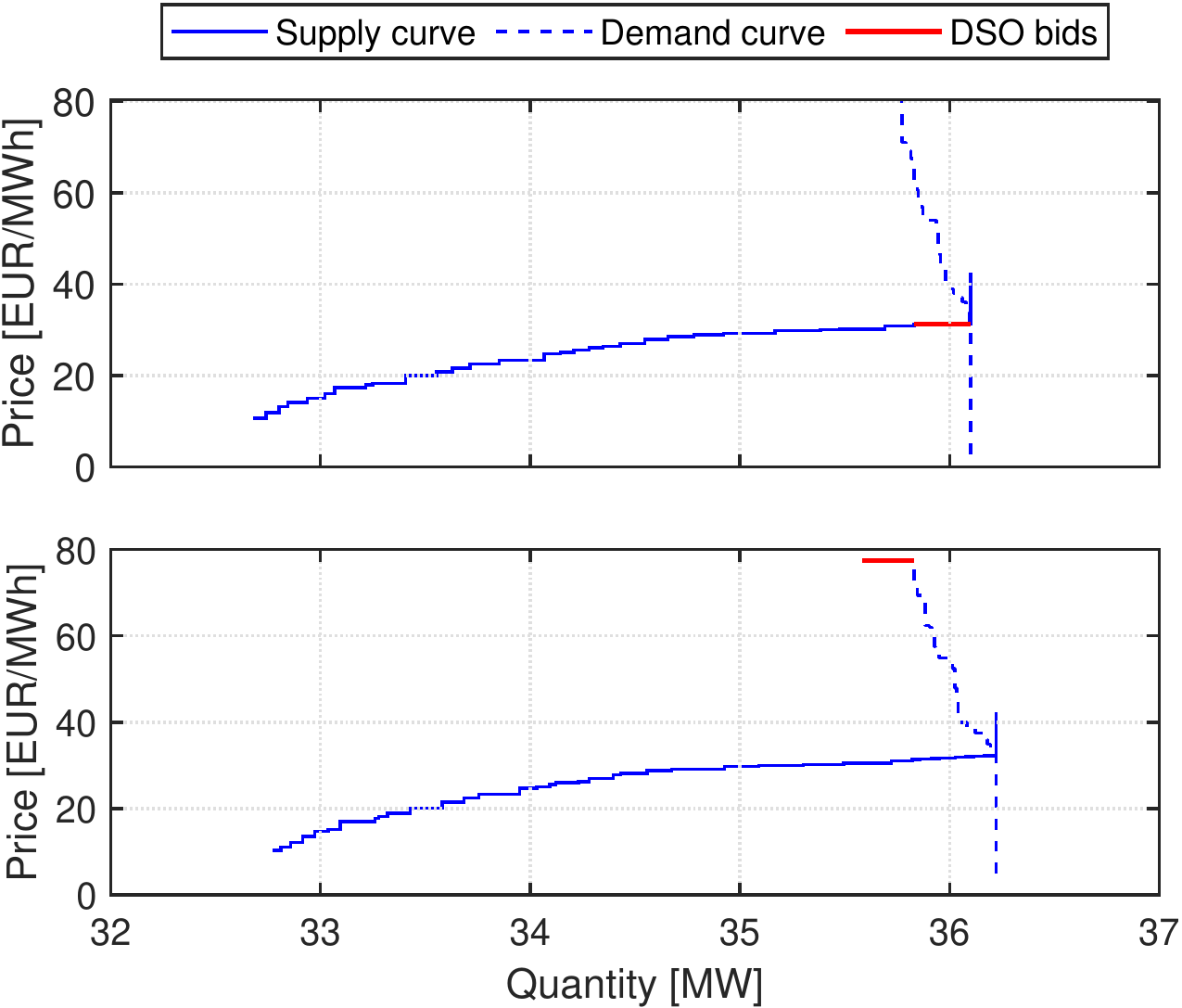}
\caption{Aggregated DA market supply curve of rival producers, demand curve, and the DA bidding/offering price-quantity pair of the DSO (i.e., aggregator) marked as red of the market for two example time slots: the top figure shows the case when the DSO acts as a net supplier and submits supply offers into the market, while the bottom figure shows the case when the DSO acts as a net consumer and submits demand bids into the market. 
        }
        \label{fig:MCcurve_example}
\end{figure}

We first validate the effectiveness of the DA model by presenting the DA bidding strategy of the DSO.
Then we conduct sensitivity analyses to investigate the impacts of the standard deviations of the day-ahead output forecast errors. Fig.~\ref{fig:MCcurve_example} shows the strategic bidding/offering strategy of the DSO for two example time slots. It can be seen that the DSO offers production at the market clearing price and bids for demand at very high market price, so as to maximize the selling price and in times of purchasing satisfying the demand of the distribution system. This is due to the assumptions that the supply of the DSO can be curtailed at cost zero while the demand is enforced to be satisfied for each time step.

\begin{figure}
     \centering
         \centering
         \includegraphics[width=0.7\textwidth]{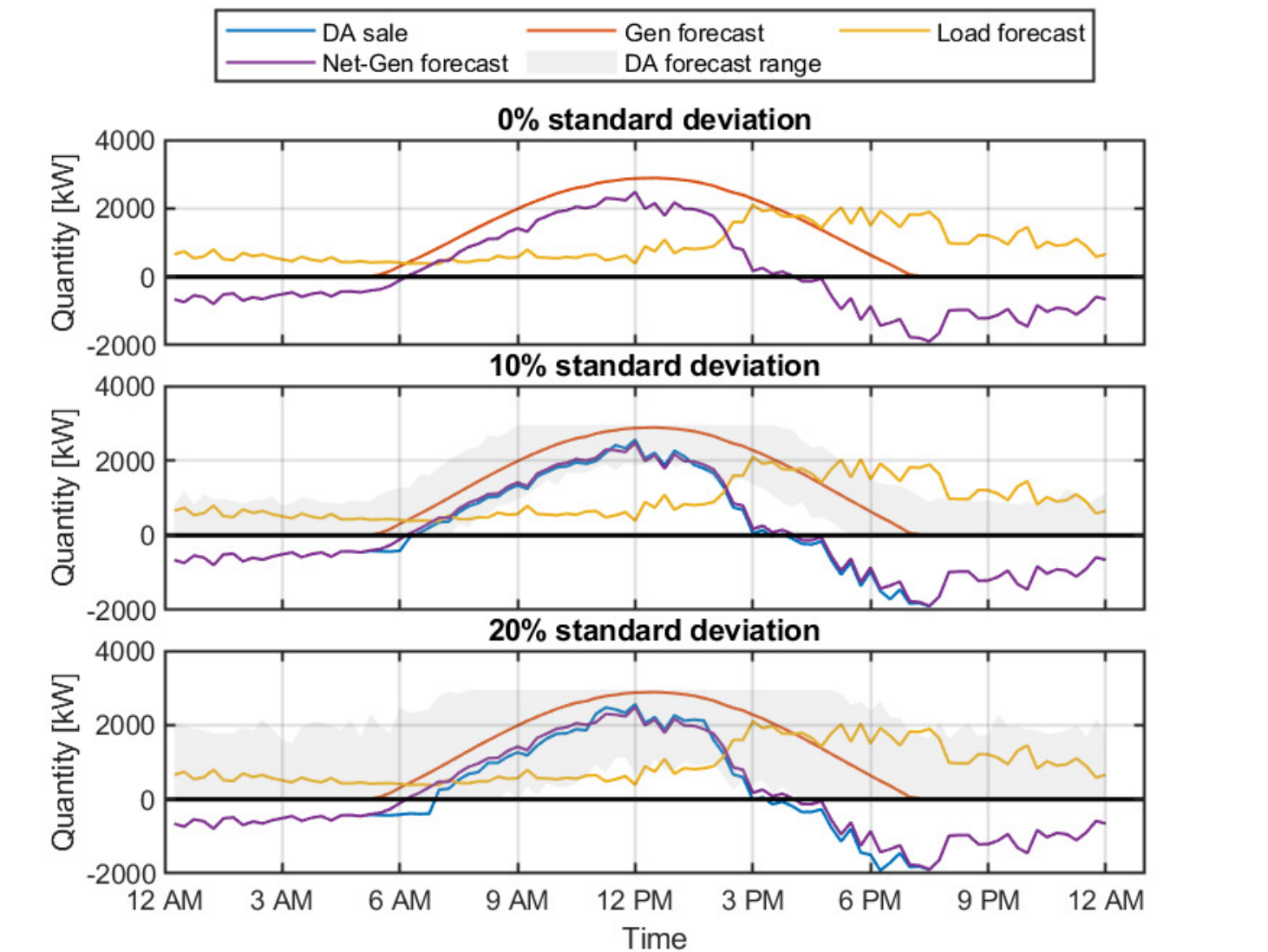}
         \caption{DA dispatch of the DSO with forecast errors generated under Gaussian distributions with mean zero and standard deviations ranging from 0 to 20\% with a step of 10\% of the generation capacity (from the top to the bottom figures), respectively.}
         \label{fig:DAdisp_SD_all}
\end{figure}

To investigate how the day-ahead bidding strategies of the DSO are impacted by the variance of the DSO's generation forecast errors, simulations are carried out using different levels of forecast errors. 
Fig. \ref{fig:DAdisp_SD_all} shows the resulting DSO's DA offering/bidding quantities setting the standard deviations of forecast errors as 0, 10\% and 20\% of the PV generation capacity. It is obvious that in general the DSO offers less to the market with the increasing standard deviation, as the worst distribution of the forecast errors worsens. However, exceptions can be observed for hours around 12 PM, when the DA market offering quantities with the standard deviation equaling 10\% or 20\% are higher than the offering quantity under zero standard deviation i.e. the forecasted net generation. This is likely due to the fact that 1) as the generation is high during these hours, with the increasing standard deviation the forecast errors after adjustment are more limited compared to the rest of the hours; 2) the DA prices of these hours are relatively higher, resulting in less penalty for negative imbalances according to equation (\ref{eq:price_imb}), i.e. a lower ratio of the negative imbalance price to the corresponding DA price compared to other hours of the day. 

The dispatched DA bidding/offering quantities with different risk aversion settings now serve as inputs to the RT market, whose results are illustrated in the following section.

\subsection{Real-time Market Results}
\begin{figure}[!tbhp]
    \centering
    \includegraphics[width=3.5in]{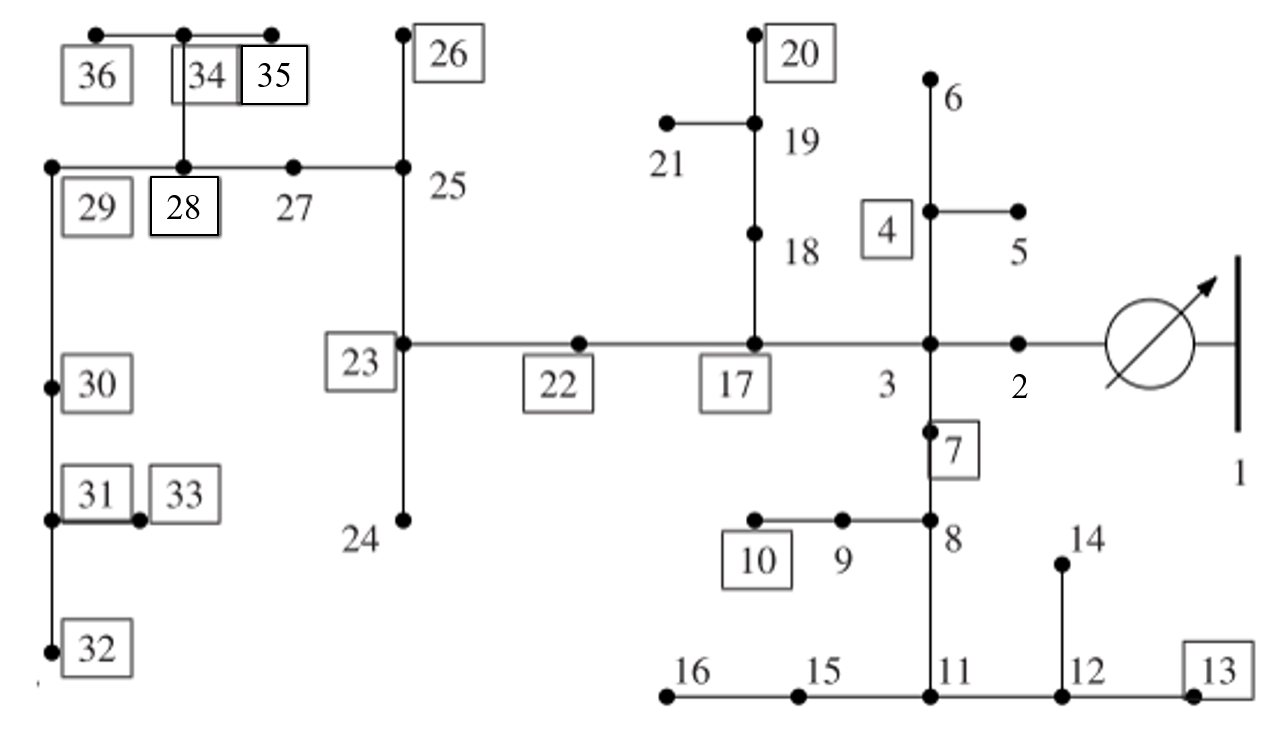}
    \caption{IEEE 37-node test feeders with 18 PV systems.}
    \label{fig:37node}
\end{figure}

After the dispatched DA supply/demand quantities have been decided in the wholesale market, we here demonstrate the proposed RT balancing market for the considered DSO of a distribution network with multiple local DERs. The considered network is a modified IEEE-37 node test feeder: the parameters of the network, such as line impedances and shunt admittances, are taken from \cite{zimmerman2010matpower}. Fig.~ \ref{fig:37node} gives the single phase equivalent of the modified network with high penetration of distributed PV systems, i.e., we place 18 PV systems in the network and their locations are marked by black boxes. Their available power is proportional to the irradiance data with 5-second granularity taken from \cite{solardata}, i.e., $\Delta T^\textrm{RT} = 5$ seconds. The original load profiles are replaced by real measurements (in 5-second resolution) from feeders in Anatolia, California, during the week of August, 2012 provided in \cite{bank2013development}. The ratings of the inverters are 200 kVA, except for the inverter at node $3$ which is 340 kVA, and at nodes $15$ and $16$ which are 200 kVA. The cost functions of DERs are defined as $J_{i,k}^{\textrm{Cost}} = a_{i,k}(p_{i,k} - p_{i,k}^{\textrm{PV}})^2 + b_{i,k}q_{i,t}^2$, which minimizes the deviation of the active power set-point $p_{i,k}$ from the PV maximum available power $p_{i,k}^{\textrm{PV}}$ and the costs of reactive power generation. The cost function parameters are set to $a_{i,k} = 3$ and $b_{i,k} = 1, \forall i\in\Omega$. The voltage limits $\underline{v}$ and $\overline{v}$ are 0.95 p.u. and 1.045 p.u..

\begin{figure}[!tbhp]
    \centering
    \includegraphics[width=5in]{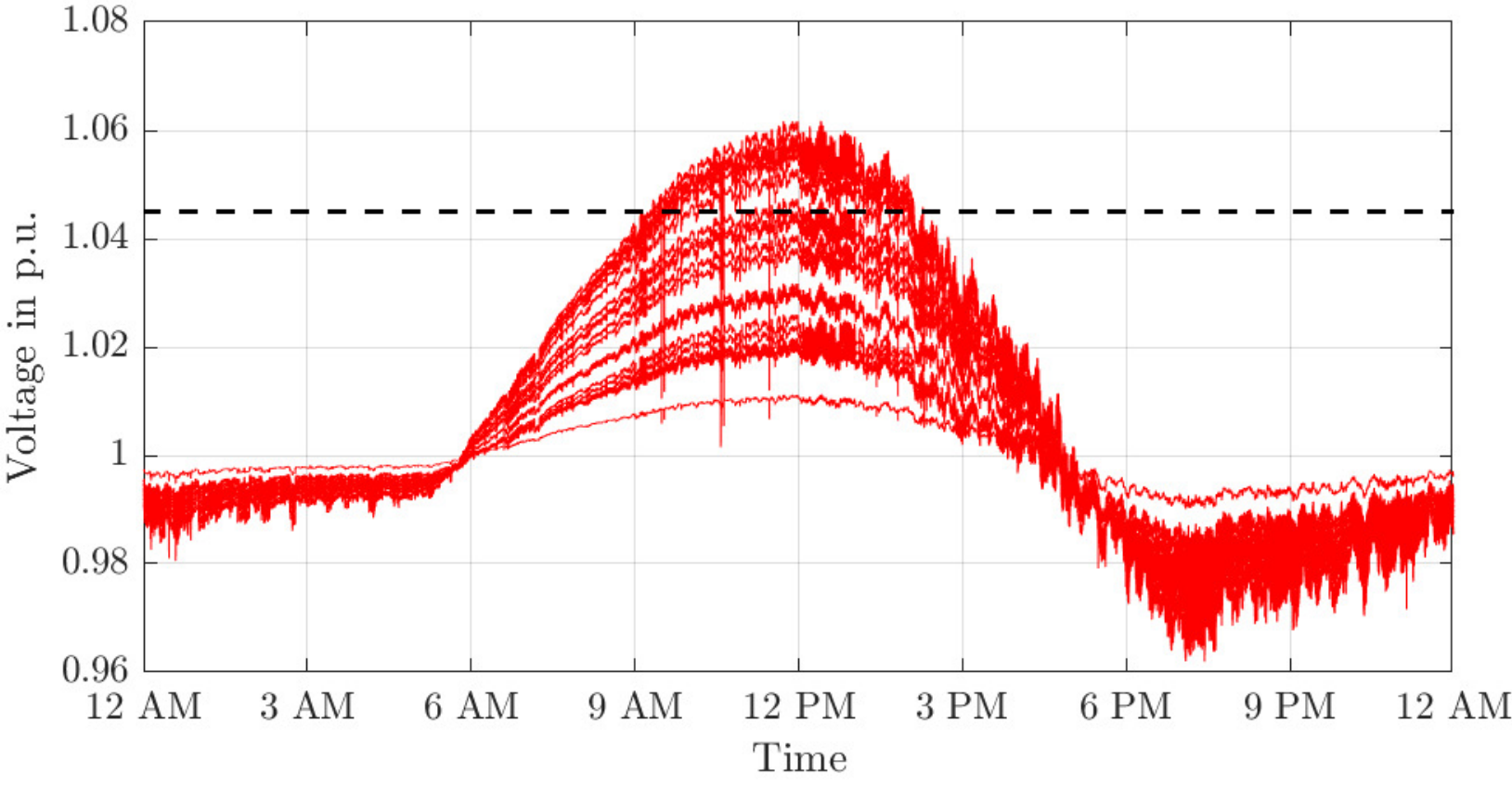}
    \caption{Voltage profile for the case without any control actions}
    \label{fig:overvoltage}
\end{figure}

\begin{figure}[!tbhp]
    \centering
    \includegraphics[width=5in]{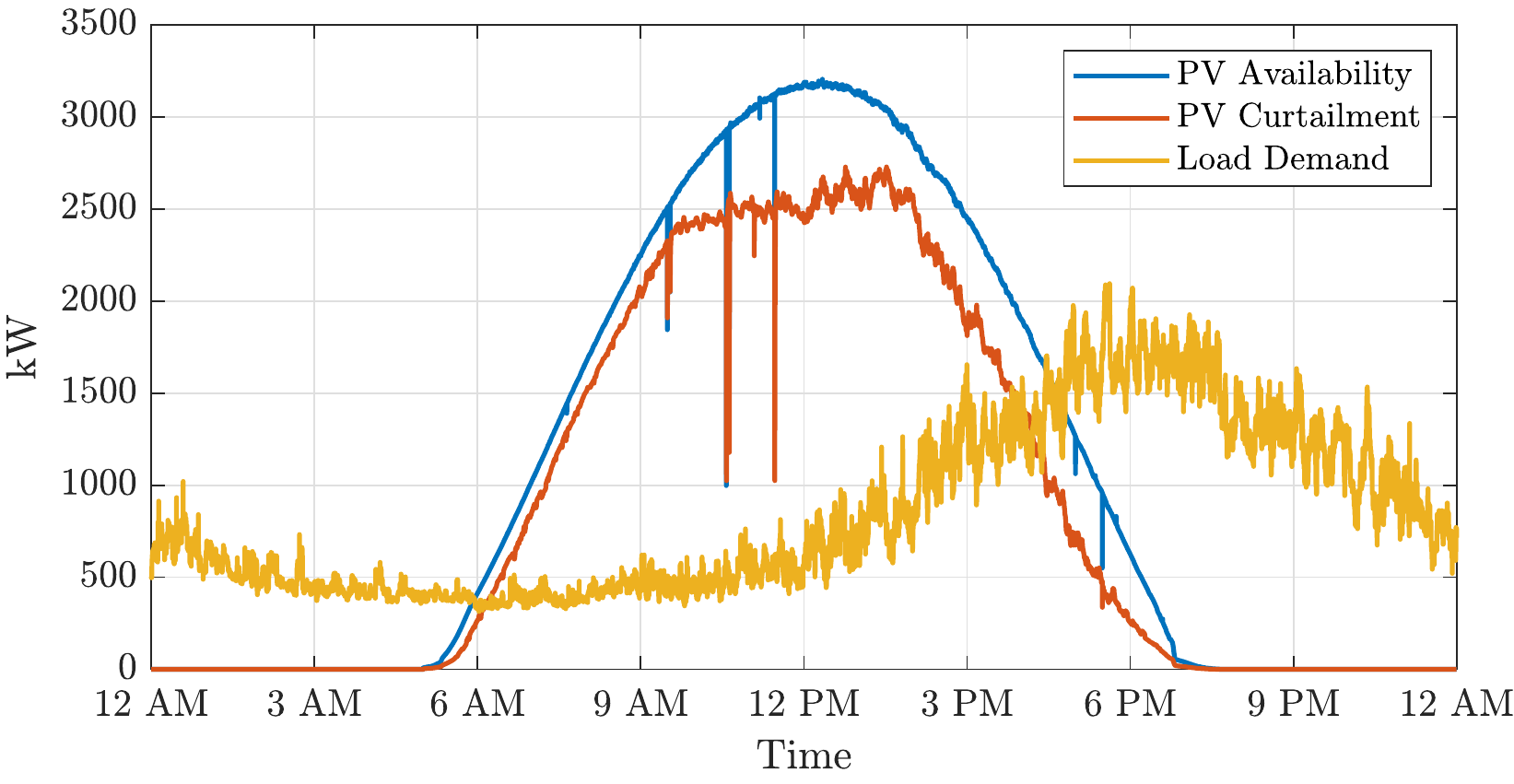}
    \caption{Aggregated solar energy availability and load demands. The solar energy curtailment is also given here after the overvoltage situation has been resolved.}
    \label{fig:pv_injections}
\end{figure}

The DA market outcomes from the previous section are used as an input to the RT market and we investigate the implications of the different standard deviations also discussed in Section \ref{subsec:DA}.  We then apply our RT market framework with the goal to minimize the deviations of the RT aggregated dispatch from the DA-ahead dispatch decisions. The default voltage of the system, i.e., without any control, is given in Fig.~\ref{fig:overvoltage}. The aggregated PV generation and loads are shown in Fig.~\ref{fig:pv_injections}. Due to high PV penetrations, an overvoltage situation can emerge during the peak of solar production. 

\begin{figure}[!tbhp]
    \centering
    \includegraphics[width=5in]{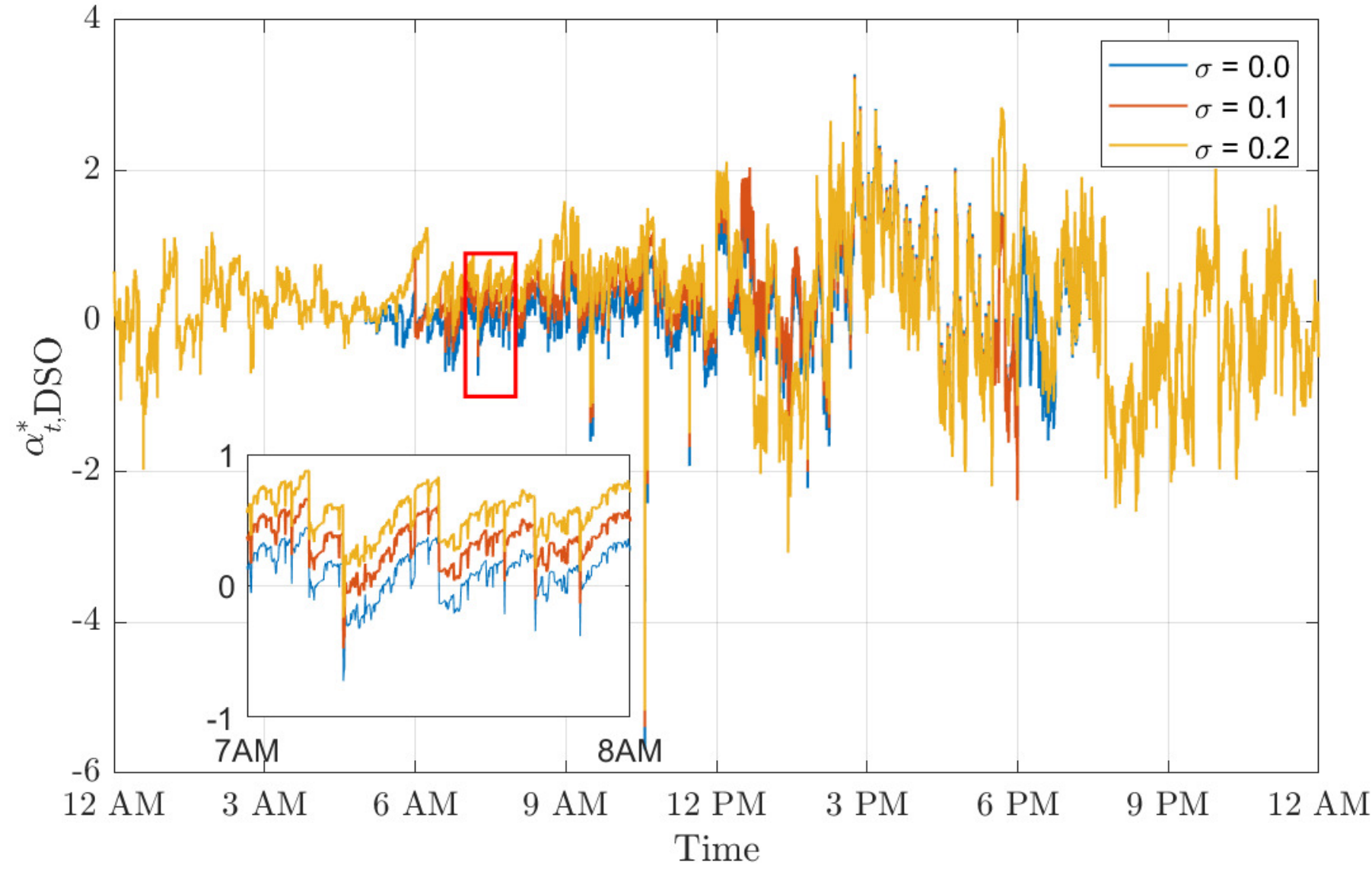}
    \caption{Comparison on incentive signal $\alpha_{\textrm{DSO},t}^*$ for the marking balancing at node 4, i.e., $\gamma = 0.2$.}
    \label{fig:alphadso}
\end{figure}

Fig.~\ref{fig:alphadso} visualizes the resulting time-varying RT incentive signals $\alpha_k^{\textrm{DSO}}$ related to the network-oriented objective for different standard deviations of the DA solar forecasts. As the standard deviation increases from $\sigma = 0$, i.e. perfect foresight, to higher values, it can be readily seen that as variations of PV forecasts $\sigma$ increase, the fluctuation of the DA purchase/selling increase as well. The online framework enables to counter larger standard deviations by generating larger incentive signals to encourage/discourage DERs' power injections. Notice that more conservative DA decisions (i.e., for larger $\sigma$) lead to more aggressive incentive signals for compensating the forecast errors of PV generations.

\begin{figure}[!tbhp]
    \centering
    \includegraphics[width=5in]{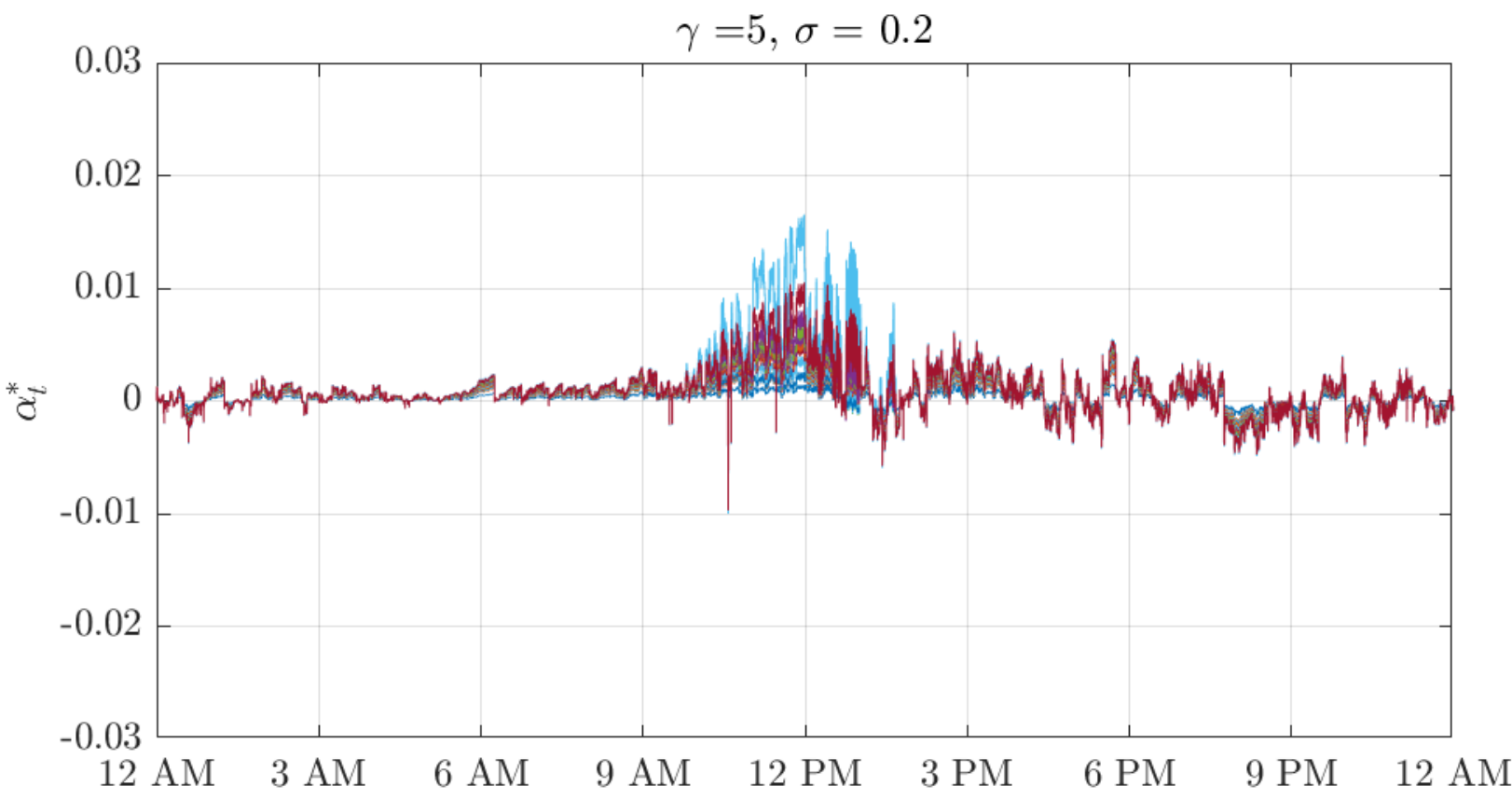}
    \caption{RT incentive signals $\alpha_{k}$ (i.e., voltage regulation \& balancing market) for all nodes with $\gamma = 5$ and $\sigma = 0.2$.}
    \label{fig:alpha_30}
\end{figure}

\begin{figure}[!tbhp]
    \centering
    \includegraphics[width=5in]{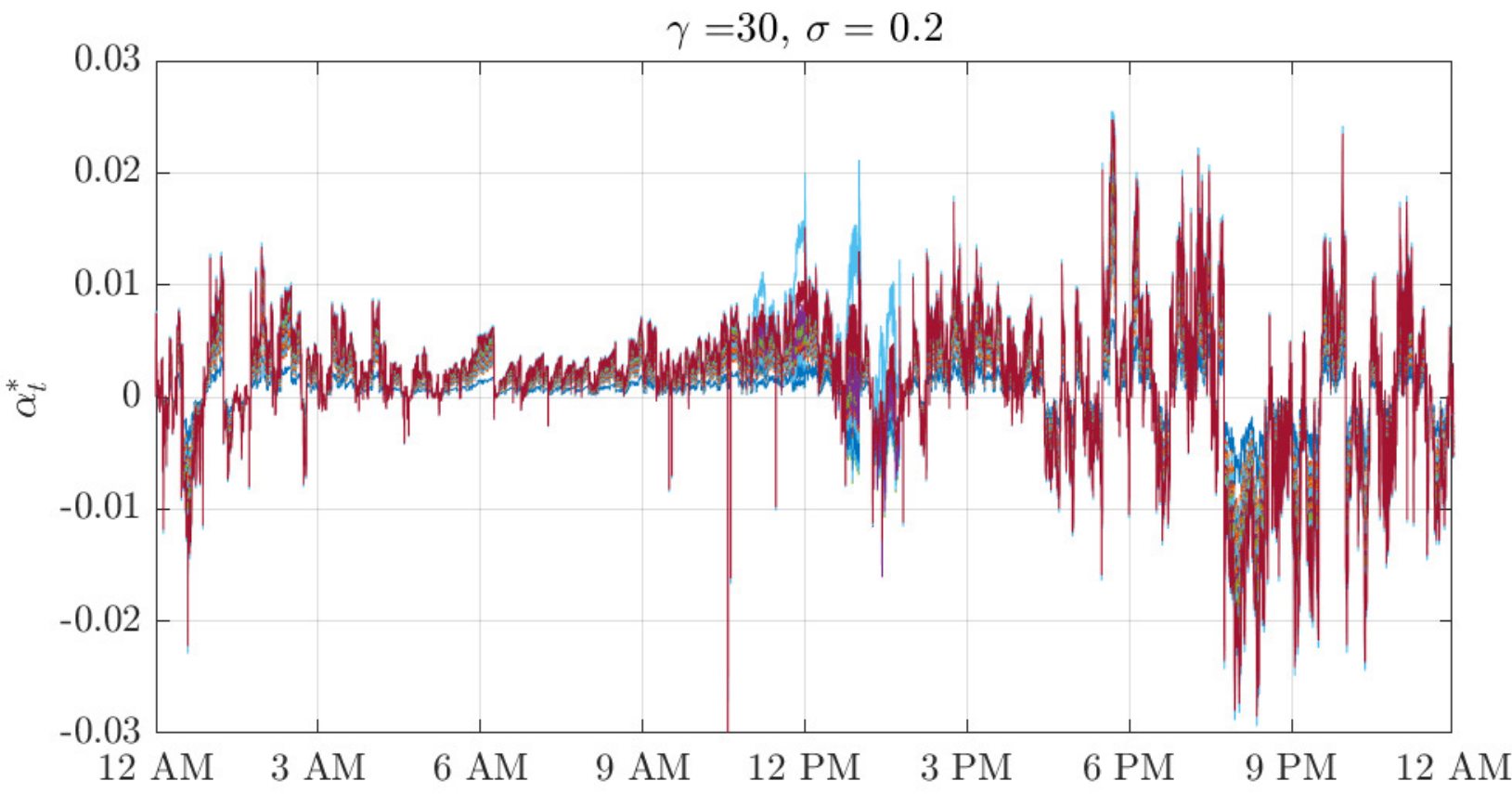}
    \caption{RT incentive signals $\alpha_{k}$ (i.e., voltage regulation \& balancing market) for all nodes with $\gamma = 30$ and $\sigma = 0.2$.}
    \label{fig:alpha_5}
\end{figure}

\begin{figure}[!tbhp]
    \centering
    \includegraphics[width=5in]{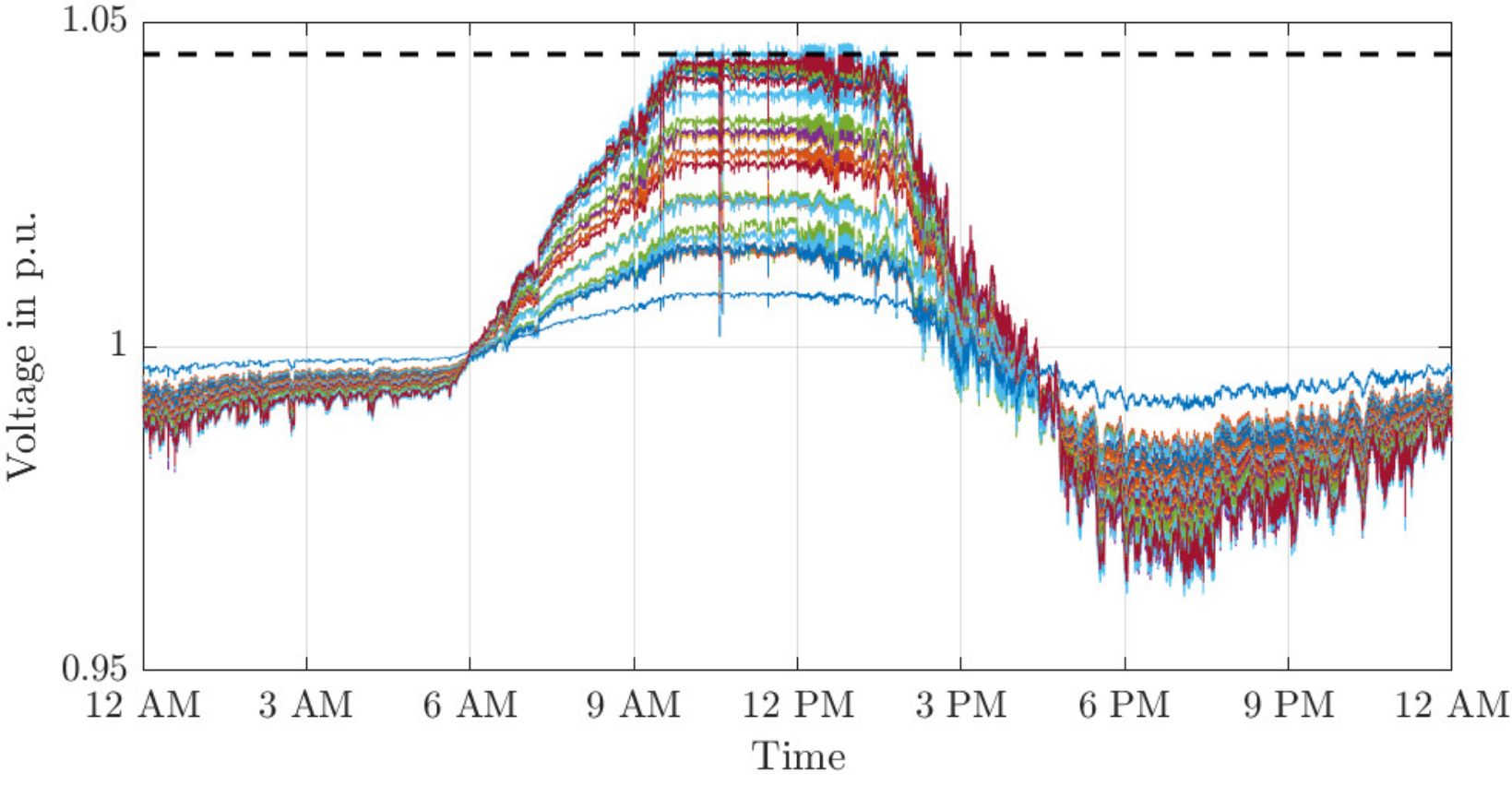}
    \caption{RT controlled voltage trajectories with $\gamma = 30$ and $\sigma = 0.2$.}
    \label{fig:voltage_with_control}
\end{figure}

Fig.~\ref{fig:alpha_30} and Fig.~\ref{fig:alpha_5} show the aggregated (for balancing and voltage regulation tasks) incentive signals $\alpha_k$ for weight factors $\gamma$ equal to 5 and 10, respectively, under the uncertainty realization, $\sigma = 0.2$. As $\gamma$ increases, the imbalance objective is emphasized. As a result, the signals for market balancing $\alpha_k^{\textrm{DSO}}$ dominates the overall incentive signal compared to the incentive signals for the voltage regulation $\alpha_k^{V}$. These parameters offer distribution system operators explicit tuning knobs to systematically design the RT market mechanism to achieve a certain network performance. Fig. \ref{fig:voltage_with_control} gives the voltages for the case of $\gamma = 30$ and $\sigma = 0.2$ indicating that voltage violations are avoided. The overvoltage has been successfully resolved under other settings of parameters, i.e., $\gamma$ and $\sigma$. In summary, we conclude that the proposed electricity market design is able to systematically consider and timely track the variations inherent to renewables in the DA wholesale market and local RT balancing market, respectively. The benefits of having a two-stage market framework can be observed in the successful trade-offs between the renewable forecast errors, network-oriented and customer-oriented objectives, while satisfying the network voltage constraints.


\section{Conclusions}\label{sec:conclusions}
\color{black}A two-stage electricity market framework is proposed in this paper. The uncertainties from DERs are handled by different optimization techniques for the different time scales. In the DA market, the bidding strategies of the aggregated DERs are optimized based on the sampled forecasting dataset via distributionally robust optimization. The controllable conservativeness of market decisions enables the TSO to operate the system taking into account different levels of risk aversion. Computational efficiency is achieved by leveraging the linear decision rule to reformulate the original bi-level problem. In the RT market, the proposed algorithm uses optimal dynamic tariffs to guide the DERs to achieve the DA decision. Note that the incentive signals show the nature of stochasticity so as to cope with large renewable variations and yet guarantee fulfilling the voltage constraints. For future work, as the generation suppliers are shifting from transmission systems to distribution networks, an interesting extension of this work is to consider the proposed market mechanism with a view on energy planning, i.e., how the proposed market mechanism could be supported by or contribute to the design of long-range market policies that foster the energy transition and help guide the regional utilities with a high penetration of DERs. \color{black}



\FloatBarrier

\bibliographystyle{elsarticle-num}

\bibliography{reference} 

\end{document}